# ON DIMENSION FOLDING OF MATRIX- OR ARRAY-VALUED STATISTICAL OBJECTS

By Bing Li[1], Min Kyung Kim and Naomi Altman

*Pennsylvania State University*

We consider dimension reduction for regression or classification in which the predictors are matrix- or array-valued. This type of predictor arises when measurements are obtained for each combination of two or more underlying variables—for example, the voltage measured at different channels and times in electroencephalography data. For these applications, it is desirable to preserve the array structure of the reduced predictor (e.g., time versus channel), but this cannot be achieved within the conventional dimension reduction formulation. In this paper, we introduce a dimension reduction method, to be called *dimension folding*, for matrix- and array-valued predictors that preserves the array structure. In an application of dimension folding to an electroencephalography data set, we correctly classify 97 out of 122 subjects as alcoholic or nonalcoholic based on their electroencephalography in a cross-validation sample.

**1. Introduction.** In many contemporary statistical applications, the sampling unit of data is in the form of a matrix- or array-valued object, such as an image, a video clip or an electroencephalography (EEG). Such data sets share two distinct characteristics: they are large, usually containing gigabytes of information, and they are structured, with each dimension of the random arrays (e.g., the rows and columns of a random matrix) representing information of a different nature. The exploration, reduction, comprehension and analysis of such large data sets, treating each array as an observation while preserving its structure, produce a fresh challenge for data analysis. In this paper, we propose a new method, to be called *dimension folding*, to deal with such types of data sets.

Received February 2009; revised July 2009.

[1]Supported in part by NSF Grants DMS-07-04621 and DMS-08-06058.

Our method is motivated by a study of an EEG data set which concerns the relationship between genetic predisposition and tendency for alcoholism ( http://kdd.ics.uci.edu/databases/eeg/eeg.data.html ). The study involved two groups of subjects: an alcoholic and a control group. Each subject was exposed to a stimulus while voltage values were measured from 64 channels of electrodes placed on the subject's scalp for 256 time points. The full data set requires about 3 gigabytes of memory. We are interested in the association between alcoholism and the pattern of voltage over times and channels.

Figure 1 shows a typical EEG pattern for an alcoholic (upper panel) and a nonalcoholic (lower panel) subject, where time and channel are represented by two horizontal axes and voltage is represented by the vertical axis. It is clear that the EEG has different patterns for the two groups. We would like to represent these different patterns in low dimension for better comprehension and classification.

In mathematical terms, the predictor is a random matrix $\mathbf{X}$ of dimension $p_L \times p_R$, and the response is a random variable $Y$—in this case, a binary random variable indicating whether or not a subject is alcoholic. We are interested in reducing the dimension of $\mathbf{X}$ as much as possible while preserving the (nonparametric) regression relation between $Y$ and $\mathbf{X}$. Without any structural restriction on the reduced predictor, the dimension reduction problem is no different from the conventional dimension reduction for vector-valued predictors. That is, one can simply treat the matrix $\mathbf{X}$ as a vector and consider the problem

$$(1) \qquad\qquad Y \perp\!\!\!\perp \text{vec}(\mathbf{X}) | \boldsymbol{\eta}^T \text{vec}(\mathbf{X}).$$

Here, $\text{vec}(\mathbf{X})$ denotes the $p_L p_R$-dimensional vector obtained by stacking the columns of $\mathbf{X}$ and $\boldsymbol{\eta}$ is a $p_L p_R \times d$ nonrandom matrix with $d < p_L p_R$. This is the classical dimension reduction problem to which all of the existing methods apply; see, for example, Li (1991, 1992), Duan and Li (1991), Cook and Weisberg (1991), Cook (1994, 1996, 1998).

However, there are practical reasons not to treat the matrix $\mathbf{X}$ as the vector $\text{vec}(\mathbf{X})$. First, problem (1) does not preserve the original matrix structure of the predictor, and so important aspects of interpretation may be lost. For example, for the EEG data, each column of $\mathbf{X}$ represents a time point and each row represents a channel. It would be desirable for the reduced predictors to still represent time and channel so that, for example, we can locate particular channels or time patterns that characterize the alcoholic tendency most distinctively. But a predictor of the form $\boldsymbol{\eta}^T \text{vec}(\mathbf{X})$ will have lost such an interpretation. Second, treating $\mathbf{X}$ as a matrix rather than a vector greatly reduces the number of parameters needed in dimension reduction, which enhances the accuracy of the estimated predictor.

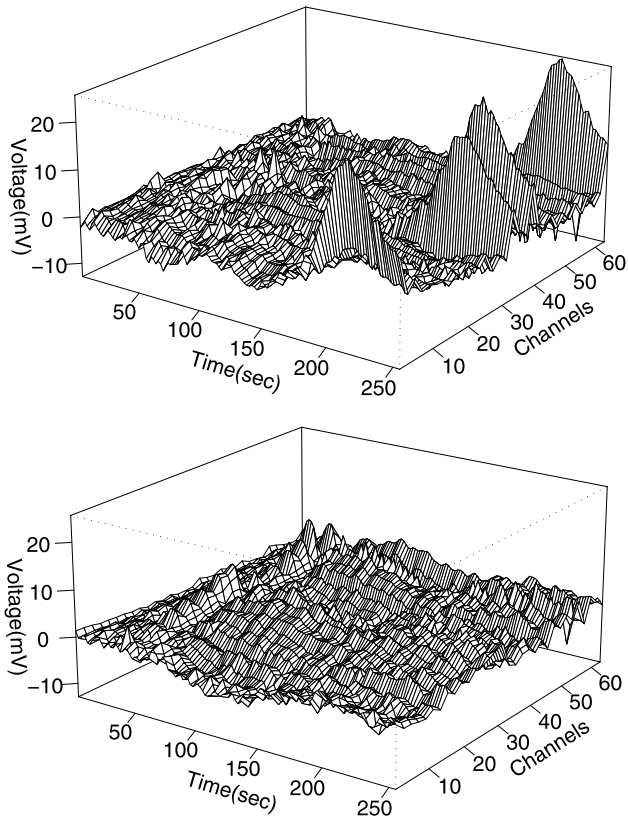

Fig. 1. *Perspective plots for the alcoholic group (upper panel) and the control group (lower panel).*

In this paper, we give a theoretical formulation and develop estimation procedures for dimension reduction problems with matrix- or array-valued predictors, which preserve the interpretations of the underlying variables. Suppose that there are matrices $\boldsymbol{\alpha}$ and $\boldsymbol{\beta}$, each with more rows than columns, such that $Y$ is independent of $\mathbf{X}$ given $\boldsymbol{\alpha}^T \mathbf{X} \boldsymbol{\beta}$. In symbols,

$$(2) \qquad Y \perp\!\!\!\perp \mathbf{X} | \boldsymbol{\alpha}^T \mathbf{X} \boldsymbol{\beta}.$$

We then only need to know the smaller matrix $\boldsymbol{\alpha}^T \mathbf{X} \boldsymbol{\beta}$ to predict, or classify, $Y$. Meanwhile, $\boldsymbol{\alpha}^T \mathbf{X} \boldsymbol{\beta}$ preserves the interpretations of channels and times—its rows representing linear combinations of channels, or principal channels, and columns representing linear combinations of times, or principal times. Such information is clearly helpful: for example, we can use the linear coefficients of the principal channel(s) to assess which parts of the brain is associated with alcoholism.

Letting $\otimes$ denote the Kronecker product, relation (2) is equivalent to

$$(3) \qquad Y \perp\!\!\!\perp \text{vec}(\mathbf{X})|(\boldsymbol{\beta} \otimes \boldsymbol{\alpha})^T \text{vec}(\mathbf{X}).$$

The challenge of dimension reduction problem (2) is that the matrix $\boldsymbol{\eta}$ in (1) cannot, in general, be written as the matrix $\boldsymbol{\beta} \otimes \boldsymbol{\alpha}$ in (3). The essence of our approach is to seek the smallest dimensional $\boldsymbol{\alpha}^T \mathbf{X} \boldsymbol{\beta}$ so that (i) $\boldsymbol{\eta}^T \mathbf{X}$ is measurable with respect to $\boldsymbol{\alpha}^T \mathbf{X} \boldsymbol{\beta}$ and (ii) the conditional independence (2) is preserved. We will also extend our results to array-valued predictors. We refer to our method as *dimension folding* to emphasize its array-preserving nature and to distinguish it from the conventional dimension reduction methods for vector-valued predictors.

In Section 2, we present the theoretical formulation and development of dimension folding. In Section 3, we introduce the key notion of the Kronecker envelope, which provides the guiding principle for constructing dimension-folding estimators from conventional dimension reduction estimators. In Sections 4 and 5, we develop three basic dimension-folding techniques: folded sliced inverse regression, folded sliced average variance estimation, and folded directional regression. In Section 6, we outline the extension to array-valued predictors. In Section 7, we make simulation comparisons between different dimension-folding methods, and between dimension-folding methods and conventional dimension reduction methods. In Section 8, we apply dimension folding to the aforementioned EEG data.

**2. Dimension-folding subspaces.** First, let us introduce some notation and terminology. For a $p \times q$ matrix $\mathbf{A}$, span($\mathbf{A}$) stands for the subspace of $\mathbb{R}^p$ spanned by the columns of $\mathbf{A}$ and $\mathbf{P_A}$ stands for the orthogonal projection onto span($\mathbf{A}$), that is, $\mathbf{P_A} = \mathbf{A}(\mathbf{A}^T\mathbf{A})^\dagger \mathbf{A}^T$, where $\dagger$ denotes the Moore–Penrose inversion. If $\mathcal{S}$ is a subspace of $\mathbb{R}^p$ and $\mathbf{A}$ is a matrix of full column rank such that span($\mathbf{A}$) $= \mathcal{S}$, then we say that $\mathbf{A}$ is a *basis matrix* of $\mathcal{S}$. Moreover, $\mathbf{P}_{\mathcal{S}}$ stands for the projection onto $\mathcal{S}$, that is, $\mathbf{P}_{\mathcal{S}} = \mathbf{P_A}$, where $\mathbf{A}$ is any basis matrix of $\mathcal{S}$. For a positive integer $p$, $\mathbf{I}_p$ denotes the $p \times p$ identity matrix.

Suppose that there are matrices $\boldsymbol{\alpha} \in \mathbb{R}^{p_L \times q_L}$ and $\boldsymbol{\beta} \in \mathbb{R}^{p_R \times q_R}$, with $q_L \leq p_L$ and $q_R \leq p_R$, such that (2) holds. This is then equivalent to

$$Y \perp\!\!\!\perp \mathbf{X}|(\boldsymbol{\alpha}\mathbf{A}_L)^T \mathbf{X}(\boldsymbol{\beta}\mathbf{A}_R),$$

whenever $\mathbf{A}_L \in \mathbb{R}^{q_L \times q_L}$ and $\mathbf{A}_R \in \mathbb{R}^{q_R \times q_R}$ are nonsingular. In other words, relation (2) depends on $\boldsymbol{\alpha}$ and $\boldsymbol{\beta}$ only through their respective column spaces, span($\boldsymbol{\alpha}$) and span($\boldsymbol{\beta}$). Thus, the identifiable parameters of this problem are column spaces of $\boldsymbol{\alpha}$ and $\boldsymbol{\beta}$, rather than $\boldsymbol{\alpha}$ and $\boldsymbol{\beta}$ themselves.

DEFINITION 1. If there exist a subspace $\mathcal{S}_L \subseteq \mathbb{R}^{p_L}$ and a subspace $\mathcal{S}_R$ of $\mathbb{R}^{d_R}$ such that

$$(4) \qquad Y \perp\!\!\!\perp \mathbf{X} | \mathbf{P}_{\mathcal{S}_L} \mathbf{X} \mathbf{P}_{\mathcal{S}_R},$$

then $\mathcal{S}_L$ is called a *left dimension-folding subspace for* $Y|\mathbf{X}$ and $\mathcal{S}_R$ is called a *right dimension-folding subspace for* $Y|\mathbf{X}$.

Under mild regularity conditions, it can be shown that if $\mathcal{S}_L$ and $\mathcal{S}'_L$ are two left dimension reduction spaces for $Y|\mathbf{X}$, then $\mathcal{S}_L \cap \mathcal{S}'_L$ is itself a left dimension reduction space. The same can be said of the right dimension reduction subspace. The situation here is similar to that in the classical setting of dimension reduction where, under very mild conditions, the intersection of two dimension reduction spaces is itself a dimension reduction space; see Cook ([1998](#)), Chiaromonte and Cook ([2001](#)) and Yin, Li and Cook ([2008](#)). Because of the similarity, we will omit the proof of this fact in the new context and take it for granted for the rest of the paper. This closure under intersection makes it possible to achieve maximal dimension folding because the intersection of all dimension-folding subspaces is itself a dimension-folding subspace. For two subspaces $\mathcal{S}_1$ and $\mathcal{S}_2$ in $\mathbb{R}^m$, let $\mathcal{S}_1 \otimes \mathcal{S}_2$ denote the linear subspace spanned by the vectors $\{\mathbf{v}_1 \otimes \mathbf{v}_2 : \mathbf{v}_1 \in \mathcal{S}_1, \mathbf{v}_2 \in \mathcal{S}_2\}$.

DEFINITION 2. Let $\mathcal{S}_{Y|\mathbf{X} \circ}$ (or $\mathcal{S}_{Y|\mathbf{X} \circ}$) be the intersection of all left (or right) dimension-folding subspaces for $Y|\mathbf{X}$. The subspace

$$\mathcal{S}_{Y|\mathbf{X} \circ} \otimes \mathcal{S}_{Y|\circ \mathbf{X}}$$

is called the *central dimension-folding subspace* and is written as $\mathcal{S}_{Y|\circ \mathbf{X} \circ}$.

Let $\boldsymbol{\beta}_L \in \mathbb{R}^{p_L \times d_L}$ be a basis matrix of $\mathcal{S}_{Y|\circ \mathbf{X}}$ and $\boldsymbol{\beta}_R \in \mathbb{R}^{p_R \times d_R}$ be a basis matrix of $\mathcal{S}_{Y|\mathbf{X} \circ}$. It is then easy to see that

$$\mathcal{S}_{Y|\circ \mathbf{X} \circ} = \mathrm{span}(\boldsymbol{\beta}_R) \otimes \mathrm{span}(\boldsymbol{\beta}_L) = \mathrm{span}(\boldsymbol{\beta}_R \otimes \boldsymbol{\beta}_L),$$

so the right-hand side is an equivalent definition of $\mathcal{S}_{Y|\circ \mathbf{X} \circ}$.

Henceforth, we no longer need to discuss any dimension-folding subspace that is not minimal in the sense of Definition [2](#), so, for brevity, we will refer to the central dimension-folding subspace simply as the *dimension-folding subspace*. Similarly, we will refer to the conventional central dimension reduction subspace (when $\mathbf{X}$ is a vector) simply as the *conventional dimension reduction subspace*. Let $\mathcal{S}_{Y|\mathrm{vec}(\mathbf{X})}$ be the conventional dimension reduction subspace of $Y$ versus the random vector $\mathrm{vec}(\mathbf{X})$. From $Y \perp\!\!\!\perp \mathrm{vec}(\mathbf{X}) | (\boldsymbol{\beta}_R \otimes \boldsymbol{\beta}_L)^T \mathrm{vec}(\mathbf{X})$, we see that

$$(5) \qquad \mathcal{S}_{Y|\mathrm{vec}(\mathbf{X})} \subseteq \mathcal{S}_{Y|\circ \mathbf{X} \circ}.$$

However, the opposite relation, $\mathcal{S}_{Y|\circ\mathbf{X}\circ} \subseteq \mathcal{S}_{Y|\operatorname{vec}(\mathbf{X})}$, does not generally hold. This means that if we do not wish to preserve the matrix structure of $\mathbf{X}$, then it is possible to further reduce the dimension of $\mathbf{X}$. However, $\mathcal{S}_{Y|\circ\mathbf{X}\circ}$ is the best that we can do if the reduced predictor is to preserve the matrix form. The following examples help to fix the idea.

EXAMPLE 1. Let $d_L = d_R = 2$ and $p_L = p_R = p$. The response $Y$ is a Bernoulli random variable with success probability equal to $\pi$; the conditional distribution of $\mathbf{X}$ given $Y$ is multivariate normal with conditional mean

$$E(\mathbf{X}|Y=0) = \mathbf{0}_{p\times p}, \qquad E(\mathbf{X}|Y=1) = \begin{pmatrix} \mu\mathbf{I}_2 & \mathbf{0}_{2\times(p-2)} \\ \mathbf{0}_{(p-2)\times 2} & \mathbf{0}_{(p-2)\times(p-2)} \end{pmatrix},$$

where $\mu \neq 0$ and $\mathbf{0}_{r\times s}$ is an $r \times s$ matrix with all of its elements equal to 0. The conditional variances are specified by

$$\operatorname{var}(X_{ij}|Y=0) = \begin{cases} \sigma^2, & (i,j) \in A, \\ 1, & (i,j) \notin A, \end{cases}$$

$$\operatorname{var}(X_{ij}|Y=1) = \begin{cases} \tau^2, & (i,j) \in A, \\ 1, & (i,j) \notin A, \end{cases}$$

where $\sigma \neq \tau$ and $A$ is the index set $\{(1,2),(2,1)\}$. We assume that $\operatorname{cov}(X_{ij}, X_{i'j'}) = 0$ whenever $(i,j) \neq (i',j')$.

Using Bayes' theorem, we can deduce that the conditional probability $P(Y=1|\mathbf{X})$ [and hence also $P(Y=0|\mathbf{X})$] is a function of $X_{11} + X_{22}$, $X_{12}^2$ and $X_{21}^2$. So, if we let $\mathbf{e}_i$ be the $p$-dimensional vector whose $i$th element is 1 and other elements are 0, then the conventional dimension reduction subspace $\mathcal{S}_{Y|\operatorname{vec}(\mathbf{X})}$ is spanned by the following three vectors in $\mathbb{R}^{p^2}$:

$$\mathbf{e}_1 \otimes \mathbf{e}_1 + \mathbf{e}_2 \otimes \mathbf{e}_2, \qquad \mathbf{e}_1 \otimes \mathbf{e}_2, \mathbf{e}_2 \otimes \mathbf{e}_1.$$

In the mean time, since the smallest submatrix of $\mathbf{X}$ that contains $X_{11} + X_{22}, X_{12}$ and $X_{21}$ is

$$\begin{pmatrix} X_{11} & X_{12} \\ X_{21} & X_{22} \end{pmatrix},$$

the central dimension-folding subspace $\mathcal{S}_{Y|\circ\mathbf{X}\circ}$ is spanned by

(6) $$\mathbf{e}_1 \otimes \mathbf{e}_1, \qquad \mathbf{e}_1 \otimes \mathbf{e}_2, \qquad \mathbf{e}_2 \otimes \mathbf{e}_1, \qquad \mathbf{e}_2 \otimes \mathbf{e}_2.$$

Thus, in this case, $\mathcal{S}_{Y|\operatorname{vec}(\mathbf{X})}$ is a proper subspace of $\mathcal{S}_{Y|\circ\mathbf{X}\circ}$.

The next example illustrates a situation where $\mathcal{S}_{Y|\operatorname{vec}(\mathbf{X})}$ and $\mathcal{S}_{Y|\circ\mathbf{X}\circ}$ coincide.

EXAMPLE 2. If we choose the index set $A$ in the definition of $\mathrm{var}(X_{ij}|Y)$ in Example 1 to be $\{(1,1),(1,2),(2,1)\}$, then it can be shown that $P(Y = 1|\mathbf{X})$ is a function of $X_{11}^2, X_{12}^2, X_{21}^2$ and $X_{11}/\tau^2 + X_{22}$. Thus, both $\mathcal{S}_{Y|\mathrm{vec}(\mathbf{X})}$ and $\mathcal{S}_{Y|\circ\mathbf{X}\circ}$ are spanned by the set of vectors in (6).

The subspace $\mathcal{S}_{Y|\circ\mathbf{X}\circ}$ enjoys an invariance property similar to that of a conventional dimension reduction subspace; see Cook (1998), Proposition 6.4.

PROPOSITION 1. *Let* $\mathbf{Z} = \mathbf{A}^T\mathbf{X}\mathbf{B}$, *where* $\mathbf{A}$ *and* $\mathbf{B}$ *are nonsingular matrices in* $\mathbb{R}^{p_L \times p_L}$ *and* $\mathbb{R}^{p_R \times p_R}$, *respectively. Then*

$$\mathcal{S}_{Y|\circ\mathbf{Z}\circ} = (\mathbf{B}^{-1} \otimes \mathbf{A}^{-1})\mathcal{S}_{Y|\circ\mathbf{X}\circ}.$$

PROOF. Let $\boldsymbol{\beta}_L$ and $\boldsymbol{\beta}_R$ be basis matrices $\mathcal{S}_{Y|\circ\mathbf{X}}$ and $\mathcal{S}_{Y|\mathbf{X}\circ}$, respectively. Because $\mathbf{Z}$ and $\mathbf{X}$ have one-to-one correspondence, we have the following equivalences:

$$Y \perp\!\!\!\perp \mathbf{X}|\boldsymbol{\beta}_L^T\mathbf{X}\boldsymbol{\beta}_R \quad \Leftrightarrow \quad Y \perp\!\!\!\perp \mathbf{X}|\boldsymbol{\beta}_L^T\mathbf{A}^{-T}\mathbf{A}^T\mathbf{X}\mathbf{B}\mathbf{B}^{-1}\boldsymbol{\beta}_R$$
$$\Leftrightarrow \quad Y \perp\!\!\!\perp \mathbf{Z}|(\mathbf{A}^{-1}\boldsymbol{\beta}_L)^T\mathbf{Z}\mathbf{B}^{-1}\boldsymbol{\beta}_R.$$

Thus, $\mathrm{span}(\mathbf{A}^{-1}\boldsymbol{\beta}_L) = \mathbf{A}^{-1}\mathcal{S}_{Y|\circ\mathbf{X}}$ is a left dimension reduction space for $Y|\mathbf{Z}$ and $\mathrm{span}(\mathbf{B}^{-1}\boldsymbol{\beta}_R) = \mathbf{B}^{-1}\mathcal{S}_{Y|\mathbf{X}\circ}$ is a right dimension reduction space for $Y|\mathbf{Z}$. Consequently,

$$\mathcal{S}_{Y|\circ\mathbf{Z}} \subseteq \mathbf{A}^{-1}\mathcal{S}_{Y|\circ\mathbf{X}}, \qquad \mathcal{S}_{Y|\mathbf{Z}\circ} \subseteq \mathbf{B}^{-1}\mathcal{S}_{Y|\mathbf{X}\circ}.$$

By the same argument, $\mathcal{S}_{Y|\circ\mathbf{X}} \subseteq \mathbf{A}\mathcal{S}_{Y|\circ\mathbf{Z}}$ and $\mathcal{S}_{Y|\mathbf{X}\circ} \subseteq \mathbf{B}\mathcal{S}_{Y|\mathbf{Z}\circ}$, which completes the proof. □

**3. Kronecker envelopes and dimension folding.** We now introduce the notion of the Kronecker envelope of a random matrix, which plays a key role in constructing dimension-folding estimators.

THEOREM 1. *Let* $\mathbf{U}$ *be an* $(r_R r_L) \times k$ *random matrix for some positive integers* $r_L$, $r_R$ *and* $k$. *There then exist subspaces* $\mathcal{S}_{\circ\mathbf{U}}$ *and* $\mathcal{S}_{\mathbf{U}\circ}$ *of* $\mathbb{R}^{r_R}$ *and* $\mathbb{R}^{r_L}$, *respectively, such that:*

1. $\mathrm{span}(\mathbf{U}) \subseteq \mathcal{S}_{\mathbf{U}\circ} \otimes \mathcal{S}_{\circ\mathbf{U}}$ *almost surely;*
2. *if there exists another pair of subspaces* $\mathcal{S}_R \in \mathbb{R}^{r_R}$ *and* $\mathcal{S}_L \in \mathbb{R}^{r_L}$ *that satisfies condition 1, then* $\mathcal{S}_{\mathbf{U}\circ} \otimes \mathcal{S}_{\circ\mathbf{U}} \subseteq \mathcal{S}_R \otimes \mathcal{S}_L$.

The random matrix $\mathbf{U}$, as well as the integers $r_L, r_R, k$, are related to specific dimension-folding methods to be described later. For example, for folded-SIR, $\mathbf{U} = \boldsymbol{\Sigma}^{-1}E[\mathrm{vec}(\mathbf{X})|Y]$, in which case $r_R = p_R$, $r_L = p_L$ and $k = 1$.

For folded-SAVE, $\mathbf{U}$ is the random matrix $\boldsymbol{\Sigma}^{-1} - \boldsymbol{\Sigma}^{-1} \operatorname{var}(\operatorname{vec}(\mathbf{X})|Y)\boldsymbol{\Sigma}^{-1}$. In this case, $r_R = p_R$, $r_L = p_L$ and $k = p_L p_R$.

PROOF OF THEOREM 1. First, we note that there always exist $\mathcal{S}_R \subseteq \mathbb{R}^{r_R}$ and $\mathcal{S}_L \subseteq \mathbb{R}^{r_L}$ so that $\operatorname{span}(\mathbf{U}) \subseteq \mathcal{S}_R \otimes \mathcal{S}_L$ because we can simply take $\mathcal{S}_R = \mathbb{R}^{r_R}$ and $\mathcal{S}_L = \mathbb{R}^{r_L}$. Thus, the following collection of subspaces is nonempty:

$$\mathfrak{F} = \{\mathcal{S}_R \otimes \mathcal{S}_L : \operatorname{span}(\mathbf{U}) \subseteq \mathcal{S}_R \otimes \mathcal{S}_L, \mathcal{S}_R \subseteq \mathbb{R}^{r_R}, \mathcal{S}_L \subseteq \mathbb{R}^{r_L}\}.$$

We will show that $\mathfrak{F}$ is a $\pi$-system [Billingsley (1986), page 36], that is, the intersection of any two members of $\mathfrak{F}$ is a member of $\mathfrak{F}$.

Let $\mathcal{S}_R \otimes \mathcal{S}_L$ and $\tilde{\mathcal{S}}_R \otimes \tilde{\mathcal{S}}_L$ be two members of $\mathfrak{F}$. Evidently, $\operatorname{span}(\mathbf{U}) \subseteq (\mathcal{S}_R \otimes \mathcal{S}_L) \cap (\tilde{\mathcal{S}}_R \otimes \tilde{\mathcal{S}}_L)$. We now show that

$$(7) \qquad (\mathcal{S}_R \otimes \mathcal{S}_L) \cap (\tilde{\mathcal{S}}_R \otimes \tilde{\mathcal{S}}_L) = (\mathcal{S}_R \cap \tilde{\mathcal{S}}_L) \otimes (\mathcal{S}_R \cap \tilde{\mathcal{S}}_L).$$

For two orthogonal subspaces, say $\mathcal{S}, \mathcal{S}'$, we use $\mathcal{S} \oplus \mathcal{S}'$ to denote the subspace spanned by the vectors in $\mathcal{S}$ and $\mathcal{S}'$. Let $\mathbf{P}_R, \tilde{\mathbf{P}}_R, \mathbf{P}_R^*$ be the projections onto $\mathcal{S}_R, \tilde{\mathcal{S}}_R, \mathcal{S}_R \cap \tilde{\mathcal{S}}_R$, respectively, and let $\mathbf{P}_L, \tilde{\mathbf{P}}_L, \mathbf{P}_L^*$ be the projections on to $\mathcal{S}_L, \tilde{\mathcal{S}}_L, \mathcal{S}_L \cap \tilde{\mathcal{S}}_L$, respectively. Then

$$\begin{aligned}
\mathcal{S}_R \otimes \mathcal{S}_L &= [\mathbf{P}_R^* \mathcal{S}_R \oplus (\mathbf{P}_R - \mathbf{P}_R^*) \mathcal{S}_R] \otimes [\mathbf{P}_L^* \mathcal{S}_L \oplus (\mathbf{P}_L - \mathbf{P}_L^*) \mathcal{S}_L] \\
&= (\mathbf{P}_R^* \mathcal{S}_R \otimes \mathbf{P}_L^* \mathcal{S}_L) \oplus [\mathbf{P}_R^* \mathcal{S}_R \otimes (\mathbf{P}_L - \mathbf{P}_L^*) \mathcal{S}_L] \\
&\quad \oplus [(\mathbf{P}_R - \mathbf{P}_R^*) \mathcal{S}_R \otimes \mathbf{P}_L^* \mathcal{S}_L] \\
&\quad \oplus [(\mathbf{P}_R - \mathbf{P}_R^*) \mathcal{S}_R \otimes (\mathbf{P}_L - \mathbf{P}_L^*) \mathcal{S}_L] \\
&\equiv (\mathbf{P}_R^* \mathcal{S}_R \otimes \mathbf{P}_L^* \mathcal{S}_L) \oplus \mathcal{A}.
\end{aligned}$$

Similarly,

$$\begin{aligned}
\tilde{\mathcal{S}}_R \otimes \tilde{\mathcal{S}}_L &= (\mathbf{P}_R^* \tilde{\mathcal{S}}_R \otimes \mathbf{P}_L^* \tilde{\mathcal{S}}_L) \oplus [\mathbf{P}_R^* \tilde{\mathcal{S}}_R \otimes (\tilde{\mathbf{P}}_L - \mathbf{P}_L^*) \tilde{\mathcal{S}}_L] \\
&\quad \oplus [(\tilde{\mathbf{P}}_R - \mathbf{P}_R^*) \tilde{\mathcal{S}}_R \otimes \mathbf{P}_L^* \tilde{\mathcal{S}}_L] \\
&\quad \oplus [(\tilde{\mathbf{P}}_R - \mathbf{P}_R^*) \tilde{\mathcal{S}}_R \otimes (\tilde{\mathbf{P}}_L - \mathbf{P}_L^*) \tilde{\mathcal{S}}_L] \\
&\equiv (\mathbf{P}_R^* \tilde{\mathcal{S}}_R \otimes \mathbf{P}_L^* \tilde{\mathcal{S}}_L) \oplus \mathcal{B}.
\end{aligned}$$

Note that

$$(8) \qquad \mathbf{P}_R^* \mathcal{S}_R \otimes \mathbf{P}_L^* \mathcal{S}_L = \mathbf{P}_R^* \tilde{\mathcal{S}}_R \otimes \mathbf{P}_L^* \tilde{\mathcal{S}}_L = (\mathcal{S}_R \cap \tilde{\mathcal{S}}_R) \otimes (\mathcal{S}_L \cap \tilde{\mathcal{S}}_L).$$

We claim that there is no nonzero common element of $\mathcal{A}$ and $\mathcal{B}$, that is, $\mathcal{A} \cap \mathcal{B} = \{\mathbf{0}\}$. This is because, by construction,

$$[(\mathbf{P}_R - \mathbf{P}_R^*) \mathcal{S}_R] \cap [(\tilde{\mathbf{P}}_R - \mathbf{P}_R^*) \mathcal{S}_R] = \{\mathbf{0}\},$$
$$[(\mathbf{P}_L - \mathbf{P}_L^*) \mathcal{S}_L] \cap [(\tilde{\mathbf{P}}_L - \mathbf{P}_L^*) \mathcal{S}_L] = \{\mathbf{0}\}.$$

It follows that

$$[(\mathbf{P}_R^* \mathcal{S}_R \otimes \mathbf{P}_L^* \mathcal{S}_L) \oplus \mathcal{A}] \cap [(\mathbf{P}_R^* \tilde{\mathcal{S}}_R \otimes \mathbf{P}_L^* \tilde{\mathcal{S}}_L) \oplus \mathcal{B}] = \mathbf{P}_R^* \mathcal{S}_R \otimes \mathbf{P}_L^* \mathcal{S}_L.$$

This, combined with (8), proves equality (7). Hence, $\mathfrak{F}$ is a $\pi$-system.

Let $\mathcal{S}_{\mathbf{U} \circ} \otimes \mathcal{S}_{\circ \mathbf{U}}$ be any member of $\mathfrak{F}$ that has the smallest dimension. It then satisfies condition 1 of the theorem. Let $\mathcal{S}_R \otimes \mathcal{S}_L$ be any member of $\mathfrak{F}$. Then $(\mathcal{S}_R \otimes \mathcal{S}_L) \cap (\mathcal{S}_{\mathbf{U} \circ} \otimes \mathcal{S}_{\circ \mathbf{U}})$ is also a member of $\mathfrak{F}$. Hence,

$$\dim[(\mathcal{S}_R \otimes \mathcal{S}_L) \cap (\mathcal{S}_{\mathbf{U} \circ} \otimes \mathcal{S}_{\circ \mathbf{U}})] = \dim(\mathcal{S}_{\mathbf{U} \circ} \otimes \mathcal{S}_{\circ \mathbf{U}}),$$

which implies that $\mathcal{S}_{\mathbf{U} \circ} \otimes \mathcal{S}_{\circ \mathbf{U}} \subseteq \mathcal{S}_R \otimes \mathcal{S}_L$. Thus, $\mathcal{S}_{\mathbf{U} \circ} \otimes \mathcal{S}_{\circ \mathbf{U}}$ satisfies condition 2, which completes the proof. $\square$

This theorem justifies the following definition of a Kronecker envelope.

DEFINITION 3. The Kronecker product space $\mathcal{S}_{\mathbf{U} \circ} \otimes \mathcal{S}_{\circ \mathbf{U}}$ in Theorem 1 is called the *Kronecker envelope of* $\mathbf{U}$ and is written as $\mathcal{E}^{\otimes}(\mathbf{U})$.

Theorem 1 guarantees that $\mathcal{E}^{\otimes}(\mathbf{U})$ exists and is uniquely defined. Note that a Kronecker envelope is defined with respect to fixed positive integers $r_L$ and $r_R$. Therefore, a fully rigorous terminology should be "Kronecker envelope of $\mathbf{U}$ with respect to integers $(r_L, r_R)$." However, in our subsequent discussions, $r_L$ and $r_R$ will be clear from the context—they will be the numbers of rows and columns of a random matrix from which $\mathbf{U}$ is derived. For this reason, we will drop this qualification.

Note that a vector $\mathbf{v} \in \mathbb{R}^{r_R r_L}$ is orthogonal to span$(\mathbf{U})$ almost surely if and only if $E[(\mathbf{v}^T \mathbf{U})^2] = \mathbf{v}^T E[(\mathbf{U} \mathbf{U}^T)]\mathbf{v} = 0$. Hence, span$[E(\mathbf{U} \mathbf{U}^T)]$ is the smallest linear subspace that contains the random subspace span$(\mathbf{U})$ almost surely. If we use $\mathcal{S}_{\mathbf{U}}$ to denote span$[E(\mathbf{U} \mathbf{U}^T)]$, then Theorem 1 and Definition 3 can both be stated with respect to $\mathcal{S}_{\mathbf{U}}$. Specifically, the condition "span$(\mathbf{U}) \subseteq \mathcal{S}_{\mathbf{U} \circ} \otimes \mathcal{S}_{\circ \mathbf{U}}$ almost surely" in Theorem 1 can be replaced by "$\mathcal{S}_{\mathbf{U}} \subseteq \mathcal{S}_{\mathbf{U} \circ} \otimes \mathcal{S}_{\circ \mathbf{U}}$" without changing the content of the theorem. In the following, we will say $\mathcal{E}^{\otimes}(\mathbf{U})$ is the Kronecker envelope of $\mathbf{U}$ or that of $\mathcal{S}_{\mathbf{U}}$ interchangeably.

In the context of conventional dimension reduction [where $\mathbf{X}$ is a vector and $\mathbf{\Sigma} = \text{var}(\mathbf{X})$], Cook, Li and Chiaromonte (2007) introduced the notion of $\mathbf{\Sigma}$-envelope as the smallest reducing subspace of $\mathbf{\Sigma}$ that contains the dimension reduction space $\mathcal{S}_{Y|\mathbf{X}}$; see also Cook, Li and Chiaromonte (2009). Their purpose was to preserve the eigenstructure of $\mathbf{\Sigma}$ so as to efficiently handle the singularity of $\mathbf{\Sigma}$. While the purpose and meaning of the Kronecker envelope differ from those of the $\mathbf{\Sigma}$-envelope, they both serve to impose extra structure on a dimension reduction (or folding) subspace, with the former imposing an eigenstructure and the latter imposing a Kronecker-product structure.

The next theorem is the theoretical basis for all of the dimension folding methods that will be described in the subsequent sections.

THEOREM 2. *Suppose that* $\mathbf{U}$ *is a random matrix in* $\mathbb{R}^{p_L p_R \times k}$ *such that* span$(\mathbf{U}) \subseteq \mathcal{S}_{Y | \mathrm{vec}(\mathbf{X})}$ *almost surely. Then*

$$\mathcal{E}^{\otimes}(\mathbf{U}) \subseteq \mathcal{S}_{Y | \circ \mathbf{X} \circ}$$

*and, consequently,* $\mathcal{S}_{\circ \mathbf{U}} \subseteq \mathcal{S}_{Y | \circ \mathbf{X}}$ *and* $\mathcal{S}_{\mathbf{U} \circ} \subseteq \mathcal{S}_{Y | \mathbf{X} \circ}$.

PROOF. By (5), span$(\mathbf{U}) \subseteq \mathcal{S}_{Y | \mathbf{X} \circ} \otimes \mathcal{S}_{Y | \circ \mathbf{X}}$ almost surely. Hence, $\mathcal{S}_{\circ \mathbf{U}} \subseteq \mathcal{S}_{Y | \circ \mathbf{X}}$ and $\mathcal{S}_{\mathbf{U} \circ} \subseteq \mathcal{S}_{Y | \mathbf{X} \circ}$. □

Theorem 2 means that if we can find a random vector or a random matrix $\mathbf{U}$ whose column space lies almost surely within the conventional dimension reduction space $\mathcal{S}_{Y | \mathrm{vec}(\mathbf{X})}$, then its Kronecker envelope is a subspace of the dimension-folding subspace $\mathcal{S}_{Y | \circ \mathbf{X} \circ}$. This is the fundamental principle by which we will construct estimates of $\mathcal{S}_{Y | \circ \mathbf{X} \circ}$. Many estimators for the conventional dimension reduction space, especially those based on conditional moments of $\mathbf{X}$ given $Y$, correspond to such random vectors or matrices. Thus, to estimate the dimension-folding subspace, all we need to do is to estimate the Kronecker envelope of the relevant random vectors or matrices which give rise to the conventional dimension reduction estimators.

We shall focus on three conventional dimension reduction estimators: SIR, SAVE and DR. In fact, using the same principle, we can develop dimension-folding methods in conjunction with all existing moment- (or conditional-moment-) based conventional methods, such as those developed in Zhu and Fang (1996), Bura and Cook (2001), Fung et al. (2002), Li (1992), Cook and Li (2002, 2004), Yin and Cook (2002), Ye and Weiss (2003), Ferre and Yao (2005) and Li, Zha and Chiaromonte (2005).

**4. Objective functions for Kronecker envelopes.** In this section, we introduce a general objective function whose minimization gives the Kronecker envelope, which will guide us in the construction of sample estimates of Kronecker envelopes.

4.1. *Conventional dimension reduction estimators.* We first review some basic facts about SIR, SAVE and DR in the conventional setting. Let $\mathbf{X}$ be a $p$-dimensional random vector and $\mathbf{\Sigma} = \mathrm{var}(\mathbf{X})$. Let $\boldsymbol{\beta}$ be a basis matrix of $\mathcal{S}_{Y | \mathbf{X}}$. SIR is based on the fact that if

(9) $$E(\mathbf{X} | \boldsymbol{\beta}^T \mathbf{X}) \text{ is linear in } \boldsymbol{\beta}^T \mathbf{X},$$

then the random vector

(10) $$\mathbf{\Sigma}^{-1} E(\mathbf{X} | Y)$$

belongs to $\mathcal{S}_{Y|\mathbf{X}}$ almost surely; see Li ([1991](#)). Let $(\tilde{\mathbf{X}}, \tilde{Y})$ be an independent copy of $(\mathbf{X}, Y)$. SAVE and DR are based on the fact that if, in addition to condition ([9](#)), we have that

(11) $$\text{var}(\mathbf{X}|\boldsymbol{\beta}^T\mathbf{X}) \text{ is a nonrandom matrix,}$$

then the column spaces for the random matrices

(12)
$$\boldsymbol{\Sigma}^{-1}[\boldsymbol{\Sigma} - \text{var}(\mathbf{X}|Y)] \qquad (\text{SAVE}),$$
$$\boldsymbol{\Sigma}^{-1}[2\boldsymbol{\Sigma} - E((\tilde{\mathbf{X}} - \mathbf{X})(\tilde{\mathbf{X}} - \mathbf{X})^T|Y, \tilde{Y})] \qquad (\text{DR})$$

are subspaces of $\mathcal{S}_{Y|\mathbf{X}}$ almost surely; see Cook and Weisberg ([1991](#)) and Li and Wang ([2007](#)).

Vectors in the central space are extracted by eigendecompositions corresponding to relations ([10](#)) and ([12](#)). For example, for SAVE, let

$$\mathbf{A} = E[\mathbf{I}_p - \text{var}(\mathbf{Z}|Y)]^2 \qquad \text{where } \mathbf{Z} = \boldsymbol{\Sigma}^{-1/2}\mathbf{X},$$

and let $\mathbf{v}_1, \ldots, \mathbf{v}_d$ be the eigenvalues of $\mathbf{A}$ corresponding to nonzero eigenvalues. Then $\{\boldsymbol{\Sigma}^{-1/2}\mathbf{v}_i : i = 1, \ldots, d\}$ spans (at least) a subspace of $\mathcal{S}_{Y|\mathbf{X}}$.

4.2. *General form of the objective function.* We now return to the matrix predictor case, where $\mathbf{X} \in \mathbb{R}^{p_L \times p_R}$. Let $\boldsymbol{\eta}$ be a basis matrix of the conventional dimension reduction subspace $\mathcal{S}_{Y|\text{vec}(\mathbf{X})}$. By the discussions in Section [4.1](#), if $E[\text{vec}(\mathbf{X})|\boldsymbol{\eta}^T\text{vec}(\mathbf{X})]$ is linear in $\boldsymbol{\eta}^T\text{vec}(\mathbf{X})$, then the random vector ([10](#)), with $\mathbf{X}$ replaced by $\text{vec}(\mathbf{X})$ and $\boldsymbol{\Sigma}$ redefined as $\text{var}[\text{vec}(\mathbf{X})]$, belongs to $\mathcal{S}_{Y|\text{vec}(\mathbf{X})}$ almost surely. If, in addition, $\text{var}[\text{vec}(\mathbf{X})|\boldsymbol{\eta}^T\text{vec}(\mathbf{X})]$ is nonrandom, then, with the same replacements, the column spaces of the random matrices in ([12](#)) are subspaces of $\mathcal{S}_{Y|\text{vec}(\mathbf{X})}$ almost surely. By Theorem [2](#), we can estimate the dimension-folding subspace $\mathcal{S}_{Y|\circ\mathbf{X}\circ}$ by targeting the Kronecker envelopes of the SIR, SAVE and DR estimators of $\mathcal{S}_{Y|\text{vec}(\mathbf{X})}$. We refer to the dimension-folding methods thus constructed as folded-SIR, folded-SAVE and folded-DR, respectively.

Again using the folded-SAVE to illustrate the idea, we minimize the objective function

$$E\|[\mathbf{I}_{p_R p_L} - \text{var}(\text{vec}(\mathbf{Z})|Y)] - \boldsymbol{\Sigma}^{1/2}(\mathbf{b} \otimes \mathbf{a})\mathbf{f}(Y)\|^2,$$

where $\text{vec}(\mathbf{Z}) = \boldsymbol{\Sigma}^{-1/2}\text{vec}(\mathbf{X})$, over matrices $\mathbf{a}$, $\mathbf{b}$ and matrix-valued functions $\mathbf{f}(\cdot)$. The matrix $\boldsymbol{\Sigma}^{1/2}$ in front of $(\mathbf{b} \otimes \mathbf{a})$ corresponds to the transformation of $\mathbf{v}_i$ to $\boldsymbol{\Sigma}^{-1/2}\mathbf{v}_i$ in a conventional procedure. The Kronecker product structure is imposed through the regression coefficient matrix $\mathbf{b} \otimes \mathbf{a}$. The next theorem shows that the solution to this minimization problem indeed gives the Kronecker envelope of $\boldsymbol{\Sigma}^{-1}[\boldsymbol{\Sigma} - \text{var}(\text{vec}(\mathbf{X}))]$, the object we desire. The theorem is stated sufficiently generally to cover all three methods.

Let $\mathbf{U}$ be a $p_R p_L \times k$ random matrix, $\boldsymbol{\alpha}_0$ and $\boldsymbol{\beta}_0$ be the basis matrices of $\mathcal{S}_{\circ \mathbf{U}}$ and $\mathcal{S}_{\mathbf{U}\circ}$, respectively, and $m_L$ and $m_R$ be the dimensions of $\mathcal{S}_{\circ \mathbf{U}}$ and $\mathcal{S}_{\mathbf{U}\circ}$, respectively. For positive integers $k_1$ and $k_2$, and a random vector $\mathbf{W}$ defined on $\Omega_{\mathbf{W}}$, let $L_2^{k_1 \times k_2}(\Omega_{\mathbf{W}})$ be the class of functions $\mathbf{f} : \Omega_{\mathbf{W}} \to \mathbb{R}^{k_1 \times k_2}$ such that $E\|\mathbf{f}(\mathbf{W})\|^2 < \infty$, where $\|\cdot\|$ is the Frobenius matrix norm.

THEOREM 3. *Suppose that the elements of $\mathbf{U}$ have finite variances and are measurable with respect to a random vector $\mathbf{W}$ and that $\mathbf{A}$ is a $p_R p_L \times p_R p_L$ nonrandom and nonsingular matrix. Let $(\mathbf{a}^*, \mathbf{b}^*, \mathbf{f}^*)$ be the minimizer of*

$$(13) \qquad E\|\mathbf{A}\mathbf{U} - \mathbf{A}(\mathbf{b} \otimes \mathbf{a})\mathbf{f}(\mathbf{W})\|^2$$

*over all $\mathbf{a} \in \mathbb{R}^{p_L \times m_L}$, $\mathbf{b} \in \mathbb{R}^{p_R \times m_R}$ and $\mathbf{f} \in L_2^{m_L m_R \times k}(\Omega_{\mathbf{W}})$. Then*

$$\mathrm{span}(\mathbf{b}^* \otimes \mathbf{a}^*) = \mathcal{E}^{\otimes}(\mathbf{U}).$$

PROOF. Since $\mathrm{span}(\boldsymbol{\beta}_0 \otimes \boldsymbol{\alpha}_0) = \mathcal{E}^{\otimes}(\mathbf{U})$ and the elements of $\mathbf{U}$ are measurable with respect to $\mathbf{W}$, there is a random matrix $\boldsymbol{\phi}(\mathbf{W}) \in L_2^{m_L m_R \times k}(\Omega_{\mathbf{W}})$ such that $\mathbf{U} = (\boldsymbol{\beta}_0 \otimes \boldsymbol{\alpha}_0)\boldsymbol{\phi}(\mathbf{W})$, which is equivalent to

$$\mathbf{A}\mathbf{U} = \mathbf{A}(\boldsymbol{\beta}_0 \otimes \boldsymbol{\alpha}_0)\boldsymbol{\phi}(\mathbf{W}).$$

Thus, (13) reaches its minimum 0 within the range of $(\mathbf{a}, \mathbf{b}, \mathbf{f})$ given in the theorem. This implies that any minimizer $(\mathbf{a}^*, \mathbf{b}^*, \mathbf{f}^*)$ of (13) must satisfy $\mathbf{A}(\mathbf{b}^* \otimes \mathbf{a}^*)\mathbf{f}^*(\mathbf{W}) = \mathbf{A}\mathbf{U}$ almost surely and, consequently,

$$(14) \qquad (\boldsymbol{\beta}_0 \otimes \boldsymbol{\alpha}_0)\boldsymbol{\phi}(\mathbf{W}) = (\mathbf{b}^* \otimes \mathbf{a}^*)f^*(\mathbf{W}) \qquad \text{almost surely.}$$

But this means that $\mathrm{span}(\mathbf{b}^* \otimes \mathbf{a}^*)$ contains $\mathbf{U}$ almost surely and has the same dimensions as $\mathcal{E}^{\otimes}(\mathbf{U})$. The theorem now follows from the uniqueness of the Kronecker envelope. $\square$

In the general objective function (13), the matrix $\mathbf{A}$ is $\boldsymbol{\Sigma}^{1/2}$ for all three dimension-folding estimators. The random element $\mathbf{W}$ is the random variable $Y$ for folded-SIR and folded-SAVE; it is the random vector $(Y, \tilde{Y})$ for folded-DR. The random element $\mathbf{U}$ is the random vector $\boldsymbol{\Sigma}^{-1} E[\mathrm{vec}(\mathbf{X})|Y]$ for folded-SIR; it is random matrix $\boldsymbol{\Sigma}^{-1}[\boldsymbol{\Sigma} - \mathrm{var}(\mathrm{vec}(\mathbf{X})|Y)]\boldsymbol{\Sigma}^{-1/2}$ for folded-SAVE; it is the random matrix

$$\boldsymbol{\Sigma}^{-1}\{2\boldsymbol{\Sigma} - E[(\mathrm{vec}(\mathbf{X}) - \mathrm{var}(\tilde{\mathbf{X}}))(\mathrm{vec}(\mathbf{X}) - \mathrm{var}(\tilde{\mathbf{X}}))^T|Y, \tilde{Y}]\}\boldsymbol{\Sigma}^{-1/2}$$

for folded-DR. Note that the $\mathbf{f}(Y)$ for folded-SIR is an $m_R m_L$-dimensional vector, whereas the $\mathbf{f}(Y)$ for folded-SAVE and $\mathbf{f}(Y, \tilde{Y})$ for folded-DR are $m_R m_L \times p_R p_L$ matrices.

The construction of the objective function in Theorem 3 expresses conditional mean in a minimization problem, which echoes the constructions used in Cook and Ni ([2005](#)), Li and Dong ([2009](#)) and Dong and Li ([2009](#)). This construction allows us to impose the Kronecker structure on minimization.

In the context of conventional dimension reduction, Li and Wang ([2007](#)) showed that both SAVE and DR are *exhaustive*, that is, the columns of the matrices in ([12](#)) do not lie in a proper subspace of $\mathcal{S}_{Y|\mathbf{X}}$. It is also known that SIR is not exhaustive when the relation between $Y$ and $\mathbf{X}$ contains a U-shaped trend. Meanwhile, even though SAVE is exhaustive at the population level, the sample estimate is often insensitive to monotone trend. Li and Wang ([2007](#)) give strong evidence that DR combines the advantages of both SIR and SAVE. A dimension-folding method inherits the exhaustive property from its conventional counterpart. Specifically, let $F_n$ be the empirical distribution based on the sample $(\mathbf{X}_1, Y_1), \ldots, (\mathbf{X}_n, Y_n)$ and let $F_0$ be the true distribution of $(\mathbf{X}, Y)$. We say that a matrix-valued statistics $\boldsymbol{\beta}(F_n)$ is an exhaustive estimator of a subspace $\mathcal{S}$ if $\operatorname{span}[\boldsymbol{\beta}(F_0)] = \mathcal{S}$.

PROPOSITION 2.    *If $\boldsymbol{\beta}(F_n)$ is an exhaustive estimator of the conventional central space $\mathcal{S}_{Y|\operatorname{vec}(\mathbf{X})}$, then the Kronecker envelope of $\operatorname{span}[\boldsymbol{\beta}(F_n)]$ is an exhaustive estimator of the dimension-folding space $\mathcal{S}_{Y|\circ\mathbf{X}\circ}$.*

PROOF.    Since $\boldsymbol{\beta}(F_n)$ is exhaustive, $\operatorname{span}[\boldsymbol{\beta}(F_0)] = \mathcal{S}_{Y|\operatorname{vec}(\mathbf{X})}$. Then $\mathcal{E} = \mathcal{E}^{\otimes}\{\operatorname{span}[\boldsymbol{\beta}(F_0)]\}$ is a Kronecker product space satisfying $Y \perp\!\!\!\perp \mathbf{X}|\mathbf{P}_{\mathcal{E}}\operatorname{vec}(\mathbf{X})$. It follows that $\mathcal{S}_{Y|\circ\mathbf{X}\circ} \subseteq \mathcal{E}$. In the mean time, since $\mathcal{S}_{Y|\circ\mathbf{X}\circ}$ is a Kronecker product space containing $\operatorname{vec}(\mathbf{X})$, we have $\mathcal{E} \subseteq \mathcal{S}_{Y|\circ\mathbf{X}\circ}$.    $\square$

**5. Estimation.**    In this section we develop an algorithm to minimize the sample version of the objective function ([13](#)). A very appealing property of the algorithm is that it can be broken down into iterations among three elementary steps, each being essentially least squares. This makes the minimization relatively fast and stable, even for a large number of parameters, which is extremely important for our applications, where the number of parameters is easily in the thousands.

5.1. *Population-level solution.*    We need the notion of a commutation matrix. If $\mathbf{A}$ is an $r_1 \times r_2$ matrix, then $\mathbf{K}_{r_1,r_2}$ is the unique matrix in $\mathbb{R}^{r_1 r_2 \times r_1 r_2}$ that transforms $\operatorname{vec}(\mathbf{A})$ to $\operatorname{vec}(\mathbf{A}^T)$: $\mathbf{K}_{r_1,r_2}\operatorname{vec}(\mathbf{A}) = \operatorname{vec}(\mathbf{A}^T)$. The explicit form and the properties of a commutation matrix can be found in Magnus and Neudecker ([1979](#)). Two properties that will be useful are: (a) if $\mathbf{A} \in \mathbb{R}^{r_1 \times r_2}$ and $\mathbf{B} = \mathbb{R}^{r_3 \times r_4}$, then

$$(15) \qquad \mathbf{A} \otimes \mathbf{B} = \mathbf{K}_{r_1,r_3}(\mathbf{B} \otimes \mathbf{A})\mathbf{K}_{r_4,r_2};$$

and (b) for any integer $r$, $\mathbf{K}_{r,1} = \mathbf{K}_{1,r} = \mathbf{I}_r$.

LEMMA 1. *Let $\mathbf{A}$ and $\mathbf{B}$ be matrices in $\mathbb{R}^{r_1 \times r_2}$ and $\mathbb{R}^{r_3 \times r_4}$, where $r_1, \ldots, r_4$ are positive integers. Then*

$$(16) \qquad \mathrm{vec}(\mathbf{A} \otimes \mathbf{B}) = \mathbf{\Pi}[\mathrm{vec}(\mathbf{A}) \otimes \mathrm{vec}(\mathbf{B})],$$

*where $\mathbf{\Pi} = \mathbf{I}_{r_2} \otimes [(\mathbf{I}_{r_4} \otimes \mathbf{K}_{r_1, r_3}) \mathbf{K}_{r_3 r_4, r_1}]$.*

PROOF. Since $\mathbf{A} \otimes \mathbf{B} = (\mathbf{a}_1 \otimes \mathbf{B}, \ldots, \mathbf{a}_{r_2} \otimes \mathbf{B})$, the vector $\mathrm{vec}(\mathbf{A} \otimes \mathbf{B})$ consists of $\mathrm{vec}(\mathbf{a}_1 \otimes \mathbf{B}), \ldots, \mathrm{vec}(\mathbf{a}_{r_2} \otimes \mathbf{B})$ stacked together vertically. By (15), $\mathbf{a}_i \otimes \mathbf{B} = \mathbf{K}_{r_1, r_3}(\mathbf{B} \otimes \mathbf{a}_i)$. Hence, $\mathrm{vec}(\mathbf{a}_i \otimes \mathbf{B}) = (\mathbf{I}_{r_4} \otimes \mathbf{K}_{r_1, r_3}) \mathrm{vec}(\mathbf{B} \otimes \mathbf{a}_i)$. But it is easy to see that $\mathrm{vec}(\mathbf{B} \otimes \mathbf{a}_i) = \mathrm{vec}(\mathbf{B}) \otimes \mathbf{a}_i$. Apply (15) again to obtain $\mathrm{vec}(\mathbf{B}) \otimes \mathbf{a}_i = \mathbf{K}_{r_3 r_4, r_1}(\mathbf{a}_i \otimes \mathrm{vec}(\mathbf{B}))$. Hence, $\mathrm{vec}(\mathbf{A} \otimes \mathbf{B})$ becomes

$$\begin{pmatrix} (\mathbf{I}_{r_4} \otimes \mathbf{K}_{r_1, r_3}) \mathbf{K}_{r_3 r_4, r_1}(\mathbf{a}_1 \otimes \mathrm{vec}(\mathbf{B})) \\ \vdots \\ (\mathbf{I}_{r_4} \otimes \mathbf{K}_{r_1, r_3}) \mathbf{K}_{r_3 r_4, r_1}(\mathbf{a}_{r_2} \otimes \mathrm{vec}(\mathbf{B})) \end{pmatrix},$$

which can be written as $\{\mathbf{I}_{r_2} \otimes [(\mathbf{I}_{r_4} \otimes \mathbf{K}_{r_1, r_3}) \mathbf{K}_{r_3 r_4, r_1}]\}[\mathrm{vec}(\mathbf{A}) \otimes \mathrm{vec}(\mathbf{B})]$, as desired. $\square$

In the following, $\mathbf{\Pi}$ is the matrix defined in Lemma 1 with $(r_1, r_2, r_3, r_4)$ taken to be $(p_R, m_R, p_L, m_L)$, that is, $\mathbf{\Pi} = \mathbf{I}_{m_R} \otimes [(\mathbf{I}_{m_L} \otimes \mathbf{K}_{p_R, p_L}) \mathbf{K}_{p_L m_L, p_R}]$.

THEOREM 4. *1. For fixed $\mathbf{f} \in L_2^{m_R m_L \times k}(\Omega_{\mathbf{W}})$, $\mathbf{a} \in \mathbb{R}^{p_L \times m_L}$, the minimizer of (13) over $\mathbf{b} \in \mathbb{R}^{p_R \times m_R}$ is $\mathbf{b} = [E(\mathbf{V}_2^T \mathbf{V}_2)]^{-1} E(\mathbf{V}_2^T \mathbf{V}_1)$, where*

$$(17) \qquad \mathbf{V}_1 = \mathrm{vec}(\mathbf{A}\mathbf{U}), \qquad \mathbf{V}_2 = (\mathbf{f}^T \otimes \mathbf{A})\mathbf{\Pi}[\mathrm{vec}(\mathbf{a}) \otimes \mathbf{I}_{p_R m_R}].$$

*2. For fixed $\mathbf{f} \in L_2^{m_R m_L \times k}(\Omega_{\mathbf{W}})$, $\mathbf{b} \in \mathbb{R}^{p_R \times m_R}$, the minimizer of (13) over $\mathbf{a} \in \mathbb{R}^{p_L \times m_L}$ is $\mathbf{a} = [E(\mathbf{V}_2^T \mathbf{V}_2)]^{-1} E(\mathbf{V}_2^T \mathbf{V}_1)$, where*

$$(18) \qquad \mathbf{V}_1 = \mathrm{vec}(\mathbf{A}\mathbf{U}), \qquad \mathbf{V}_2 = (\mathbf{f}^T \otimes \mathbf{A})\mathbf{\Pi}[\mathbf{I}_{p_L m_L} \otimes \mathrm{vec}(\mathbf{b})].$$

*3. For fixed $\mathbf{a} \in \mathbb{R}^{p_L \times m_L}$ and $\mathbf{b} \in \mathbb{R}^{p_R \times m_R}$, the minimizer of (13) over $\mathbf{f} \in L_2^{m_R m_L \times k}(\Omega_{\mathbf{W}})$ is $\mathbf{f}(\mathbf{w}) = (\mathbf{V}_2^T \mathbf{V}_2)^{-1}[\mathbf{V}_2^T \mathbf{V}_1(\mathbf{w})]$, where*

$$(19) \qquad \mathbf{V}_1(\mathbf{w}) = \mathrm{vec}[\mathbf{A}\mathbf{U}(\mathbf{w})], \qquad \mathbf{V}_2 = \mathbf{I}_{p_R p_L} \otimes [\mathbf{A}(\mathbf{b} \otimes \mathbf{a})].$$

PROOF. By standard calculations, if $\mathbf{V}_1$ is an $r_1$-dimensional random vector and $\mathbf{V}_2$ is an $r_1 \times r_2$-dimensional random matrix, each having finite second moments, then the minimizer of

$$(20) \qquad E\|\mathbf{V}_1 - \mathbf{V}_2 \mathbf{c}\|^2$$

over all $\mathbf{c} \in \mathbb{R}^{r_2}$ is

$$(21) \qquad \mathbf{c}^* = [E(\mathbf{V}_2^T \mathbf{V}_2)]^{-1} E(\mathbf{V}_2^T \mathbf{V}_1).$$

We now rewrite the objective function (13) as

(22) $$E\|\text{vec}(\mathbf{AU}) - \text{vec}[\mathbf{A}(\mathbf{b} \otimes \mathbf{a})\mathbf{f}]\|^2.$$

To prove part 1, note that

$$\text{vec}[\mathbf{A}(\mathbf{b} \otimes \mathbf{a})\mathbf{f}] = (\mathbf{f}^T \otimes \mathbf{A})\text{vec}(\mathbf{b} \otimes \mathbf{a}) = (\mathbf{f}^T \otimes \mathbf{A})\mathbf{\Pi}[\text{vec}(\mathbf{b}) \otimes \text{vec}(\mathbf{a})],$$

where the second equality follows from Lemma 1. Note that $\text{vec}(\mathbf{b}) \otimes \text{vec}(\mathbf{a}) = \text{vec}[\text{vec}(\mathbf{b})\,\text{vec}^T(\mathbf{a})] = [\text{vec}(\mathbf{a}) \otimes \mathbf{I}_{p_R m_R}]\text{vec}(\mathbf{b})$. Hence,

$$\text{vec}[\mathbf{A}(\mathbf{b} \otimes \mathbf{a})\mathbf{f}] = (\mathbf{f}^T \otimes \mathbf{A})\mathbf{\Pi}[\text{vec}(\mathbf{a}) \otimes \mathbf{I}_{p_R m_R}]\text{vec}(\mathbf{b}).$$

Thus, (22) is of the form (20), with $\mathbf{V}_1, \mathbf{V}_2$ defined as in (17) and $\mathbf{c} = \text{vec}(\mathbf{b})$. The assertion of part 1 now follows from (21).

The proof of part 2 is similar to that of part 1 and is thus omitted. Let us turn to part 3. For each fixed $\mathbf{w}$, $\mathbf{f}(\mathbf{w})$ is the minimizer of

(23) $$\begin{aligned} &E[\|\mathbf{AU} - \mathbf{A}(\mathbf{b} \otimes \mathbf{a})\mathbf{f}(\mathbf{Z})\|^2 | \mathbf{W} = \mathbf{w}] \\ &= \|\text{vec}[\mathbf{AU}(\mathbf{w})] - [\mathbf{I}_{p_R p_L} \otimes \mathbf{A}(\mathbf{b} \otimes \mathbf{a})]\text{vec}[\mathbf{f}(\mathbf{w})]\|^2, \end{aligned}$$

where, since $\mathbf{U}(\mathbf{w})$ and $\mathbf{f}(\mathbf{w})$ are fixed given $\mathbf{W} = \mathbf{w}$, the conditional expectation $E(\cdot | \mathbf{W} = \mathbf{w})$ disappears. Now, apply (21) to (23) with $\mathbf{V}_1(\mathbf{w}), \mathbf{V}_2(\mathbf{w})$ defined in (19) and $\mathbf{c} = \text{vec}[\mathbf{f}(\mathbf{w})]$ to complete the proof. □

When $k = 1$, as is the case for folded-SIR, the solution can be further simplified. Let $\mathbf{a}$ be a vector in $\mathbb{R}^{rs}$, where $r$ and $s$ are positive integers. Thus, $\mathbf{a}$ can be written as $(\mathbf{a}_1^T, \ldots, \mathbf{a}_s^T)^T$, where each $\mathbf{a}_i$ is a vector in $\mathbb{R}^r$. We define $\text{mat}_r(\mathbf{a})$ to be the $r \times s$ matrix $(\mathbf{a}_1, \ldots, \mathbf{a}_s)$. This is an inverse operation of vec, in the sense that, for any matrix $\mathbf{A} \in \mathbb{R}^{r \times s}$ and any vector $\mathbf{a} \in \mathbb{R}^{rs}$, we have

$$\text{mat}_r[\text{vec}(\mathbf{A})] = \mathbf{A}, \qquad \text{vec}[\text{mat}_r(\mathbf{a})] = \mathbf{a}.$$

Note that the operation $\text{mat}_r$ is specified by a number $r$, but no such specification is needed for the definition of vec. A useful property of the mat operation is that if $\mathbf{A} \in \mathbb{R}^{r_1 \times r_2}$, $\mathbf{b} \in \mathbb{R}^{r_2 r_3}$ and $\mathbf{C} \in \mathbb{R}^{r_3 \times r_4}$ for some positive integers $r_1, \ldots, r_4$, then

(24) $$\text{mat}_{r_1}[(\mathbf{C}^T \otimes \mathbf{A})\mathbf{b}] = \mathbf{A}\,\text{mat}_{r_2}(\mathbf{b})\mathbf{C}.$$

This can be verified by taking vec on both sides and observing that

$$\text{vec}[\mathbf{A}\,\text{mat}_{r_2}(\mathbf{b})\mathbf{C}] = (\mathbf{C}^T \otimes \mathbf{A})\text{vec}[\text{mat}_{r_2}(\mathbf{b})] = (\mathbf{C}^T \otimes \mathbf{A})\mathbf{b}.$$

If $k = 1$, then $\mathbf{f}$ is an $m_R m_L$-dimensional vector. So, by (24), $(\mathbf{b} \otimes \mathbf{a})\mathbf{f}$ can be written as $\text{vec}[\mathbf{a}\,\text{mat}_{m_L}(\mathbf{f})\mathbf{b}^T]$, which, in turn, can be written as

$$[\mathbf{I}_{p_R} \otimes \mathbf{a}\,\text{mat}_{m_L}(\mathbf{f})]\text{vec}(\mathbf{b}^T) = [\mathbf{I}_{p_R} \otimes \mathbf{a}\,\text{mat}_{m_L}(\mathbf{f})]\mathbf{K}_{p_R, m_R}\text{vec}(\mathbf{b})$$

or

$$[\mathbf{b} \operatorname{mat}_{m_L}^T(\mathbf{f}) \otimes \mathbf{I}_{p_L}] \operatorname{vec}(\mathbf{a}).$$

Thus, the $\mathbf{V}_2$ in (17), (18), (19) in Theorem 4 can be replaced by

$$\mathbf{A}[\mathbf{I}_{p_R} \otimes \mathbf{a} \operatorname{mat}_{m_L}(\mathbf{f})]\mathbf{K}_{p_R, m_R}, \mathbf{A}[\mathbf{b} \operatorname{mat}_{m_L}^T(\mathbf{f}) \otimes \mathbf{I}_{p_L}], \mathbf{A}(\mathbf{b} \otimes \mathbf{a}),$$

respectively. This alternative expression often requires less computation for folded-SIR.

5.2. *Numerical procedures.* We now describe the estimation procedures for folded-SIR and folded-DR at the sample level. The procedure for folded-SAVE is similar to folded-DR and is thus omitted. Let $(\mathbf{X}_1, Y_1), \ldots, (\mathbf{X}_n, Y_n)$ be an i.i.d. sample of $(\mathbf{X}, Y)$. We estimate $\mathbf{\Sigma}$ by the sample moment

$$\hat{\mathbf{\Sigma}} = n^{-1} \sum_{i=1}^{n} \operatorname{vec}(\mathbf{X}_i - \bar{\mathbf{X}}) \operatorname{vec}^T(\mathbf{X}_i - \bar{\mathbf{X}}).$$

As with the conventional dimension reduction methods such as SIR, we discretize the response $Y$. Let $J_1, \ldots, J_s$ be a partition of $\Omega_Y$. Let $D = \delta(Y)$ be the discrete random variable defined by

$$\delta(Y) = \ell \qquad \text{if } Y \in J_\ell, \ell = 1, \ldots, s.$$

For a function $h$ of $(\mathbf{X}, Y)$, let $E_n h(\mathbf{X}, Y)$ denote the sample average $n^{-1} \times \sum_{i=1}^{n} h(\mathbf{X}_i, Y_i)$. We summarize the estimating procedure for the folded-SIR as the following five-step algorithm:

1. generate the initial values of $\mathbf{a}_0 \in \mathbb{R}^{p_L \times m_L}$, $\{\mathbf{f}_0(\ell) : \ell = 1, \ldots, s\}$, say, from a sample of the $N(0, 1)$ variables;
2. for $\ell = 1, \ldots, s$, compute $\hat{p}_\ell = E_n[I(D = \ell)]$ and

$$\hat{\mathbf{V}}_1(\ell) = \hat{p}_\ell^{-1} \operatorname{vec}\{\hat{\mathbf{\Sigma}}^{-1/2} E_n[\operatorname{vec}(\mathbf{X})I(D = \ell)]\},$$

$$\hat{\mathbf{V}}_2(\ell) = \hat{\mathbf{\Sigma}}^{1/2}[\mathbf{I}_{p_R} \otimes \mathbf{a}_0 \operatorname{mat}_{m_L}(\mathbf{f}_0(\ell))]\mathbf{K}_{p_R, m_R},$$

   then compute $\operatorname{vec}(\mathbf{b}_1)$ by

(25) $$\left[\sum_{i=1}^{s} \hat{p}_\ell \hat{\mathbf{V}}_2^T(\ell)\hat{\mathbf{V}}_2(\ell)\right]^{-1}\left[\sum_{i=1}^{s} \hat{p}_\ell \hat{\mathbf{V}}_2^T(\ell)\hat{\mathbf{V}}_1(\ell)\right];$$

3. recompute $\hat{\mathbf{V}}_2(\ell)$ as $\hat{\mathbf{\Sigma}}^{1/2}[\mathbf{b}_1 \operatorname{mat}_{m_L}^T(\mathbf{f}_0(\ell)) \otimes \mathbf{I}_{p_L}]$, then compute $\operatorname{vec}(\mathbf{a}_1)$ using (25), but with the recomputed $\hat{\mathbf{V}}_2(\ell)$;
4. recompute $\hat{\mathbf{V}}_2$ as $\mathbf{I}_{p_R p_L} \otimes [\hat{\mathbf{\Sigma}}^{1/2}(\mathbf{b}_1 \otimes \mathbf{a}_1)]$, noting that, at this step, $\hat{\mathbf{V}}_2$ does not depend on $\ell$, then compute $\mathbf{f}_1(\ell) = (\hat{\mathbf{V}}_2^T \hat{\mathbf{V}}_2)^{-1}\hat{\mathbf{V}}_2^T \hat{\mathbf{V}}_1(\ell)$;

5. return to step 2 and iterate, each time using the most recent $\mathbf{a}$, $\mathbf{b}$ and $\mathbf{f}$, until

$$\sum_{\ell=1}^{s} \hat{p}_\ell \|\hat{\boldsymbol{\Sigma}}^{-1/2} E_n[\text{vec}(\mathbf{X})|D=\ell] - \hat{\boldsymbol{\Sigma}}^{1/2}(\mathbf{b} \otimes \mathbf{a})\mathbf{f}(\ell)\|^2$$

stabilizes, then use the resulting $\mathbf{a}$ and $\mathbf{b}$ as the estimates of $\boldsymbol{\alpha}_0$ and $\boldsymbol{\beta}_0$.

The algorithm for folded-DR is similar to folded-SIR, except that the $\hat{\mathbf{V}}_1$ and $\hat{\mathbf{V}}_2$ become more complicated. Let $\tilde{\boldsymbol{\nabla}} = \text{vec}(\tilde{\mathbf{X}})$, $\boldsymbol{\Delta} = \tilde{\boldsymbol{\nabla}} - \boldsymbol{\nabla}$ and $\tilde{D} = \delta(\tilde{Y})$. Then, for $k, \ell = 1, \ldots, s$,

$$E(\boldsymbol{\Delta}\boldsymbol{\Delta}^T|D=k, \tilde{D}=\ell)$$
$$= E(\boldsymbol{\nabla}\boldsymbol{\nabla}^T|D=k) - E(\boldsymbol{\nabla}|D=k)E(\tilde{\boldsymbol{\nabla}}^T|\tilde{D}=\ell)$$
$$- E(\tilde{\boldsymbol{\nabla}}|\tilde{D}=\ell)E(\boldsymbol{\nabla}^T|D=k) + E(\tilde{\boldsymbol{\nabla}}\tilde{\boldsymbol{\nabla}}^T|\tilde{D}=\ell).$$

Let $n_1, \ldots, n_s$ be the numbers of observations in slices $J_1, \ldots, J_s$. The sample estimate for the above conditional expectation is

$$E_n(\boldsymbol{\Delta}\boldsymbol{\Delta}^T|D=k, \tilde{D}=\ell)$$
$$= \frac{1}{n_k}\sum_{r \in J_k} \boldsymbol{\nabla}_r \boldsymbol{\nabla}_r^T - \frac{1}{n_k n_\ell}\sum_{r \in J_k} \boldsymbol{\nabla}_r \sum_{t \in J_\ell} \boldsymbol{\nabla}_t^T$$
$$- \frac{1}{n_k n_\ell}\sum_{r \in J_\ell} \boldsymbol{\nabla}_r \sum_{t \in J_k} \boldsymbol{\nabla}_t^T + \frac{1}{n_\ell}\sum_{t \in J_\ell} \boldsymbol{\nabla}_t \boldsymbol{\nabla}_t^T.$$

We now summarize the algorithm for folded-DR:

1. generate the initial values of $\mathbf{a}_0 \in \mathbb{R}^{p_L \times m_L}$, $\{\mathbf{f}_0(k, \ell) : k, \ell = 1, \ldots, s\}$ from, say, a sample of the $N(0, 1)$ variables;

2. for $k, \ell = 1, \ldots, s$, compute $\hat{p}_{k\ell} = n_k n_\ell / n$ and

$$\hat{\mathbf{V}}_1(k, \ell) = \text{vec}\{\hat{\boldsymbol{\Sigma}}^{-1/2}[2\hat{\boldsymbol{\Sigma}} - E_n(\boldsymbol{\Delta}\boldsymbol{\Delta}^T|D=k, D=\ell)]\hat{\boldsymbol{\Sigma}}^{-1/2}\},$$
$$\hat{\mathbf{V}}_2(k, \ell) = [\mathbf{f}_0^T(k, \ell) \otimes \boldsymbol{\Sigma}^{1/2}]\boldsymbol{\Pi}[\text{vec}(\mathbf{a}_0) \otimes \mathbf{I}_{p_R m_R}],$$

then compute $\text{vec}(\mathbf{b}_1)$ using the formula

$$(26) \quad \left[\sum_{k=1}^{s}\sum_{\ell=1}^{s} \hat{p}_{k,\ell} \hat{\mathbf{V}}_2^T(k, \ell)\hat{\mathbf{V}}_2(k, \ell)\right]^{-1}\left[\sum_{k=1}^{s}\sum_{\ell=1}^{s} \hat{p}_{k,\ell} \hat{\mathbf{V}}_2^T(k, \ell)\hat{\mathbf{V}}_1(k, \ell)\right];$$

3. recompute $\hat{\mathbf{V}}_2(k, \ell)$ as

$$\hat{\mathbf{V}}_2(k, \ell) = [\mathbf{f}_0^T(k, \ell) \otimes \boldsymbol{\Sigma}^{1/2}]\boldsymbol{\Pi}[\mathbf{I}_{p_L m_L} \otimes \text{vec}(\mathbf{b}_1)],$$

then compute $\text{vec}(\mathbf{a}_1)$ by (26), using the newly computed $\hat{\mathbf{V}}_2(k, \ell)$;

4. compute $\mathbf{f}_1(k, \ell)$ by

$$\mathbf{f}_1(k, \ell) = [(\mathbf{b}_0 \otimes \mathbf{a}_1)^T \hat{\boldsymbol{\Sigma}} (\mathbf{b}_1 \otimes \mathbf{a}_1)]^{-1} (\mathbf{b}_1 \otimes \mathbf{a}_1)^T$$
$$\times [2\hat{\boldsymbol{\Sigma}} - E_n(\boldsymbol{\Delta}\boldsymbol{\Delta}^T | D = k, \tilde{D} = \ell)] \hat{\boldsymbol{\Sigma}}^{-1/2};$$

5. repeat steps 2, 3, 4, using the most updated $\mathbf{a}$, $\mathbf{b}$ and $\mathbf{f}$ at each step, until the objective function

$$\sum_{k, \ell = 1}^s \hat{p}_{k\ell} \| \hat{\boldsymbol{\Sigma}}^{-1/2} [2\hat{\boldsymbol{\Sigma}} - E_n(\boldsymbol{\nabla}\boldsymbol{\nabla}^T | D = k, \tilde{D} = \ell)] \hat{\boldsymbol{\Sigma}}^{-1/2} - \hat{\boldsymbol{\Sigma}}^{1/2} (\mathbf{b} \otimes \mathbf{a}) \mathbf{f}(l, \tilde{\ell}) \|^2$$

stabilizes.

5.3. *Singularity of $\hat{\boldsymbol{\Sigma}}$.* When $p_R p_L > n$, the sample covariance matrix $\hat{\boldsymbol{\Sigma}}$ of vec($\mathbf{X}$) is singular and $\hat{\boldsymbol{\Sigma}}^{-1}$ does not exist. There are several ways to deal with this. One is to replace vec($\mathbf{X}$) by its principal components. Chiromonte and Martinell (2002) and Li and Li (2004) used this method in the conventional setting. If all principal components corresponding to nonzero eigenvalues of $\hat{\boldsymbol{\Sigma}}$ are used, then this amounts to using the Moore–Penrose inverse $\hat{\boldsymbol{\Sigma}}^\dagger$ in place of $\hat{\boldsymbol{\Sigma}}^{-1}$. Another option is to use the ridge-regression-type inverse $(\hat{\boldsymbol{\Sigma}} + \varepsilon \mathbf{I}_{p_R p_L})^{-1}$, where $\varepsilon > 0$, in place of $\hat{\boldsymbol{\Sigma}}^{-1}$; see Hoerl (1962) and Marquardt (1970). For a related development in conventional dimension reduction, see Tyekucheva and Chiaromonte (2008) and Li (2008). Finally, it is possible to adapt the iterative transformation approach of Cook, Li and Chiaromonte (2007) to dimension folding, but further research is needed in this regard. In the subsequent simulations and application, we use the first two approaches to handle the singularity of $\hat{\boldsymbol{\Sigma}}$.

5.4. *Robustness.* The dimension-folding methods proposed here are based on sample moments, which are known to be sensitive to outliers. Zhou (2009) described a weighting scheme to achieve robustness in the conventional setting for a dimension reduction method derived from canonical correlations [Fung et al. (2002)]. We outline how that scheme can be adapted to dimension folding.

For a given sample, let $\hat{w}(\mathbf{x})$, $\mathbf{x} \in \mathbb{R}^{p_L \times p_R}$ be a decreasing and nonnegative function of $[\text{vec}(\mathbf{x}) - E_n \text{vec}(\mathbf{x})]^T \hat{\boldsymbol{\Sigma}}^{-1} [\text{vec}(\mathbf{x}) - E_n \text{vec}(\mathbf{x})]$ such that $\sum_{i=1}^n \hat{w}(\mathbf{X}_i) = 1$. To downplay observations lying far away from the center of observed predictors, we replace the empirical measure that assigns probability mass $1/n$ to each pair $(\mathbf{X}_i, Y_i)$ with the alternative random measure that assigns probability mass $\hat{w}(\mathbf{X}_i)$ to $(\mathbf{X}_i, Y_i)$. We then replace the usual sample moments and sample conditional moments by moments calculated from this alternative measure. For example, for folded-SIR, we replace $E_n \operatorname{vec}(\mathbf{X})$, $\hat{\boldsymbol{\Sigma}}$ and $E_n[\operatorname{vec}(\mathbf{X})|D = \ell]$ by

$$E_n^* \operatorname{vec}(\mathbf{X}_i) = \sum_{i=1}^{n} \hat{w}(\mathbf{X}_i) \operatorname{vec}(\mathbf{X}_i),$$

$$\sum_{i=1}^{n} \hat{w}(\mathbf{X}_i)[\operatorname{vec}(\mathbf{X}_i) - E_n^* \operatorname{vec}(\mathbf{X})][\operatorname{vec}(\mathbf{X}_i) - E_n^* \operatorname{vec}(\mathbf{X})]^T,$$

$$\sum_{i=1}^{n} \hat{w}(\mathbf{X}_i) \operatorname{vec}(\mathbf{X}_i) I(D_i = \ell) \bigg/ \sum_{i=1}^{n} \hat{w}(\mathbf{X}_i) I(D_i = \ell).$$

The rest of the algorithm remains the same. Folded-SAVE and folded-DR can be robustified in a similar fashion.

**6. Array-valued predictors.** As mentioned in the Introduction, we sometimes also encounter sampling units in the form of higher-dimensional arrays. For example, a video clip is a three-dimensional array. In this section, we extend dimension folding to these cases. For reasons of brevity, we omit the details of algorithms, which can be constructed analogously.

Let $\mathbf{X} = \{X_{j_1 \cdots j_u} : j_1 = 1, \ldots, p_1, \ldots, j_u = 1, \ldots, p_u\}$ be a $u$-way random array of dimension $p_1 \times \cdots \times p_u$, and let $Y$ be a scalar-valued random response. Our goal is to reduce $\mathbf{X}$ to a smaller $u$-way array of dimension $d_1 \times \cdots \times d_u$ while preserving the regression relation between $\mathbf{X}$ and $Y$. That is, we seek nonrandom matrices

$$\boldsymbol{\alpha}^{(1)} = \{\alpha_{i_1 j_1}^{(1)} : i_1 = 1, \ldots, p_1, j_1 = 1, \ldots, d_1\},$$

$$\vdots$$

$$\boldsymbol{\alpha}^{(u)} = \{\alpha_{i_u j_u}^{(u)} : i_u = 1, \ldots, p_u, j_u = 1, \ldots, d_u\},$$

such that $Y$ is conditional independent of $\mathbf{X}$ given the array

$$(27) \quad \left\{ \sum_{i_1=1}^{p_1} \cdots \sum_{i_u=1}^{p_u} \alpha_{i_1 j_1}^{(1)} \cdots \alpha_{i_u j_u}^{(u)} X_{i_1 \cdots i_u} : i_1 = 1, \ldots, d_1, \ldots, i_u = 1, \ldots, d_u \right\}.$$

Let $\operatorname{vec}(\mathbf{X})$ denote the vector of elements of $\mathbf{X}$ with its first index changing the fastest. That is,

$$\operatorname{vec}(\mathbf{X}) = (X^{1 \cdots 1}, X^{2 \cdots 1}, \ldots, X^{1 \cdots 2}, X^{2 \cdots 2}, \ldots, X^{p_1 \cdots p_u})^T.$$

Parallel to the definition of the $\operatorname{mat}_r$ operator introduced in Section 5, we define $\operatorname{arr}_{p_1 \cdots p_u}$ as the inverse operator of $\operatorname{vec}(\mathbf{X})$. That is, $\operatorname{arr}_{p_1 \cdots p_u}[\operatorname{vec}(\mathbf{X})] = \mathbf{X}$. The array (27) can then be written as

$$\operatorname{arr}_{p_1 \cdots p_u}[(\boldsymbol{\alpha}^{(u)} \otimes \cdots \otimes \boldsymbol{\alpha}^{(1)})^T \operatorname{vec}(\mathbf{X})].$$

Since the above array has a one-to-one relation with $(\boldsymbol{\alpha}^{(u)} \otimes \cdots \otimes \boldsymbol{\alpha}^{(1)})^T \operatorname{vec}(\mathbf{X})$, the general dimension-folding problem can be stated as

$$(28) \qquad Y \perp\!\!\!\perp \mathbf{X} | (\boldsymbol{\alpha}^{(u)} \otimes \cdots \otimes \boldsymbol{\alpha}^{(1)})^T \operatorname{vec}(\mathbf{X}).$$

We note that, as in the matrix-predictor case, the order of $1, \ldots, u$ is reversed in the string of Kronecker products: the coefficient matrix associated with the last index of $\mathbf{X}$ appears first in the string of Kronecker products.

The central dimension-folding subspace is then defined as the smallest subspace

$$\operatorname{span}(\boldsymbol{\alpha}^{(u)}) \otimes \cdots \otimes \operatorname{span}(\boldsymbol{\alpha}^{(1)})$$

for which the relation (28) is satisfied. This subspace will be written as $\mathcal{S}_{Y|\mathbf{X}^{\circ u}}$.

Once again, the idea is to start with a random matrix whose column space lies almost surely in the conventional dimension reduction space $\mathcal{S}_{Y|\operatorname{vec}(\mathbf{X})}$ and to use the Kronecker envelope of this random matrix to estimate the dimension-folding subspace. The next theorem is parallel to Theorem 1 and Definition 3, so its proof is omitted.

DEFINITION 4.   Let $\mathbf{U}$ be a random matrix in $\mathbb{R}^{(p_1 \cdots p_u) \times k}$. There are subspaces $\mathcal{S}_1 \subseteq \mathbb{R}^{p_1}, \ldots, \mathcal{S}_u \subseteq \mathbb{R}^{p_n}$ of dimensions $t_1 \leq p_1, \ldots, t_u \leq p_u$ such that:

1. $\operatorname{span}(\mathbf{U}) \subseteq \mathcal{S}_1 \otimes \cdots \otimes \mathcal{S}_u$ almost surely;
2. if there exists another $u$-tuple of subspaces $\mathcal{S}'_1 \subseteq \mathbb{R}^{p_1}, \ldots, \mathcal{S}'_u \subseteq \mathbb{R}^{p_u}$ that satisfies condition 1, then

$$\mathcal{S}_1 \otimes \cdots \otimes \mathcal{S}_u \subseteq \mathcal{S}'_1 \otimes \cdots \otimes \mathcal{S}'_u$$

and the subspace $\mathcal{S}_1 \otimes \cdots \otimes \mathcal{S}_u$ is called the *Kronecker envelope of* $\mathbf{U}$.

We denote the generalized Kronecker envelope of $\mathbf{U}$ by $\mathcal{E}^{\otimes}_{p_1, \ldots, p_r}(\mathbf{U})$. Using an argument similar to that used in Section 3, we can prove the following result.

THEOREM 5.   *Let* $\mathbf{X}$ *be a random array in* $\mathbb{R}^{p_1 \times \cdots \times p_u}$. *If* $\mathbf{U}$ *is a random matrix in* $\mathbb{R}^{(p_1 \cdots p_u) \times k}$ *whose column space is contained in* $\mathcal{S}_{Y|\operatorname{vec}(\mathbf{X})}$ *almost surely, then* $\mathcal{E}^{\otimes}_{p_1, \ldots, p_u}(\mathbf{U})$ *is contained in* $\mathcal{S}_{Y|\mathbf{X}^{\circ u}}$.

This theorem provides the guiding principle for estimating the central dimension-folding space $\mathcal{S}_{Y|\mathbf{X}^{\circ u}}$. That is, we start with a conventional dimension reduction method such as SIR, SAVE or DR for estimating $\mathcal{S}_{Y|\operatorname{vec}(\mathbf{X})}$ and then construct the estimates of its Kronecker envelope via objective functions analogous to (13).

**7. Simulation studies.** In this section, we evaluate by simulation the performance of the three dimension folding estimators, in a classification problem in which the response is a binary variable and the predictor matrices corresponding to the two values of $Y$ differ both in location and variation. Our comparison is twofold: we compare the performance among the three dimension folding methods themselves and compare them with conventional dimension reduction methods when the dimension reduction subspace coincides with the dimension-folding subspace. While dimension-folding methods are introduced primarily to preserve the array structure, intuitively, they should be more accurate than their conventional counterparts when $\mathcal{S}_{\circ X \circ} = \mathcal{S}_{Y \mid \text{vec}(\mathbf{X})}$ because the former contains far fewer parameters. The second comparison is made in order to confirm this intuition.

To assess the accuracy of a dimension-folding method, we use the criterion

$$\|\mathbf{P}_{\hat{\boldsymbol{\beta}} \otimes \hat{\boldsymbol{\alpha}}} - \mathbf{P}_{\boldsymbol{\beta} \otimes \boldsymbol{\alpha}}\|, \tag{29}$$

where $\|\cdot\|$ is a matrix norm, which, for example, can be the Frobenius norm or the largest singular value. This is a measure of discrepancy between the subspaces $\text{span}(\boldsymbol{\beta}^* \otimes \boldsymbol{\alpha}^*)$ and $\text{span}(\boldsymbol{\beta} \otimes \boldsymbol{\alpha})$; see Li, Zha and Chiaromonte (2005) for intuition about and further discussion of this criterion. In the following, we use the Frobenius norm.

To make a sensible comparison, it is helpful to define a "benchmark" of this discrepancy, that is, its value when the two spaces are not related at all. Let $\boldsymbol{\alpha}^* \in \mathbb{R}^{p_L \times d_L}$ and $\boldsymbol{\beta}^* \in \mathbb{R}^{p_R \times d_R}$ be random matrices whose entries are i.i.d. standard normal. We define $E(\|\mathbf{P}_{\boldsymbol{\beta}^* \otimes \boldsymbol{\alpha}^*} - \mathbf{P}_{\boldsymbol{\beta} \otimes \boldsymbol{\alpha}}\|)$ to be the *benchmark distance*. The benchmark is easily computed by simulation. It depends on dimensions $p_L, p_R, d_L$ and $d_R$, but is independent of the model and the estimator, as well as of $\boldsymbol{\alpha}$ and $\boldsymbol{\beta}$ (despite its appearance). A similar benchmark was used in Li, Wen and Zhu (2008) in the classical setting. The performance of the conventional dimension reduction methods is assessed similarly. Let $\hat{\boldsymbol{\eta}}$ be the conventional dimension reduction estimates of $\boldsymbol{\eta}$, a basis matrix for $\mathcal{S}_{Y \mid \text{vec}(\mathbf{X})}$. We use

$$\|\mathbf{P}_{\hat{\boldsymbol{\eta}}} - \mathbf{P}_{\boldsymbol{\eta}}\| \tag{30}$$

to assess the error of the conventional methods. Note that $\mathbf{P}_{\boldsymbol{\beta} \otimes \boldsymbol{\alpha}} = \mathbf{P}_{\boldsymbol{\eta}}$, so the comparison is on equal footing.

EXAMPLE 1 (Continued). Let $\mathbf{X}$ and $Y$ be defined as in Example 1 in Section 2. We take $\pi = 1/2$, $\sigma^2 = 0.1$ and $\tau^2 = 1.5$. Recall that, in this case, $\mathcal{S}_{Y \mid \text{vec}(\mathbf{X})}$ is a proper subset of $\mathcal{S}_{Y \mid \circ X \circ}$.

We generate $n$ pairs of observations, $(\mathbf{X}_1, Y_1), \ldots, (\mathbf{X}_n, Y_n)$, from this model, with $n = 100, 200, 300, 500, 800$ and $p = 5, 10$. We apply folded-SIR, folded-SAVE and folded-DR. Table 1 gives the means of criterion (29), as calculated from $N = 500$ simulated samples for each combination of $n$ and $p$. The standard errors of these means are all within 0.02 and are not presented. From the table, we can see that the overall best performer is folded-DR, followed by folded-SAVE and folded-SIR. Both folded-DR and folded-SAVE perform much better than folded-SIR. This is because the two mixing components for $Y = 0$ and $Y = 1$ differ both by location and variance, the latter of which cannot be captured by folded-SIR.

TABLE 1
*Comparison among dimension-folding methods*

| Method | $n = 100$ | $n = 200$ | $n = 300$ | $n = 500$ | $n = 800$ |
|---|---|---|---|---|---|
| | $p_L = p_R = 5$ (benchmark distance = 2.586) | | | | |
| Folded-SIR | 1.115 | 0.751 | 0.598 | 0.496 | 0.369 |
| Folded-SAVE | 0.566 | 0.295 | 0.220 | 0.161 | 0.121 |
| Folded-DR | 0.531 | 0.287 | 0.215 | 0.158 | 0.119 |
| | $p_L = p_R = 10$ (benchmark distance = 2.772) | | | | |
| Folded-SIR | 2.006 | 1.289 | 1.034 | 0.772 | 0.604 |
| Folded-SAVE | 2.710 | 1.410 | 0.581 | 0.345 | 0.236 |
| Folded-DR | 2.296 | 1.019 | 0.542 | 0.331 | 0.230 |

To give a sense of how the methods perform, in Figure 2, we present the scatterplot matrices of the four elements of $\hat{\mathbf{a}}^T \mathbf{X} \hat{\mathbf{b}}$ ($m_L = m_R = 2$), as estimated by folded-SIR (left panel) and folded-DR (right panel). We see that the four predictors by folded-SIR separate the two groups by location, whereas folded-DR separates them by both location and variation, as shown in the $(X11, X22)$ plot in the lower panel.

EXAMPLE 2 (Continued). Let $\mathbf{X}$ and $Y$ be defined as in Example 2 in Section 2. Again take $\pi = 1/2$, $\sigma^2 = 0.1$ and $\tau^2 = 1.5$. The difference from the previous example is that the index set $A$ for the definition of $\text{var}(X_{ij}|Y)$ is changed to ensure that $\mathcal{S}_{Y|\text{vec}(\mathbf{X})} = \mathcal{S}_{Y|\circ\mathbf{X}\circ}$, so that the comparison of dimension-folding methods and conventional dimension reduction methods is on an equal footing.

With the same choices of $n$, $p$ and $N$, in Table 2, we present the means of either criterion (29) (for dimension-folding methods) or criterion (30) (for conventional dimension reduction methods). We observe very substantial improvements by dimension folding. This is due to the fact that the column space of $\boldsymbol{\beta} \otimes \boldsymbol{\alpha}$ contains far fewer parameters than the column space of $\boldsymbol{\eta}$, if both matrices have the same dimension. We also see that folded-SAVE and folded-DR perform much better than folded-SIR and the same pattern holds for their conventional counterparts, for the same reason explained in the previous comparison.

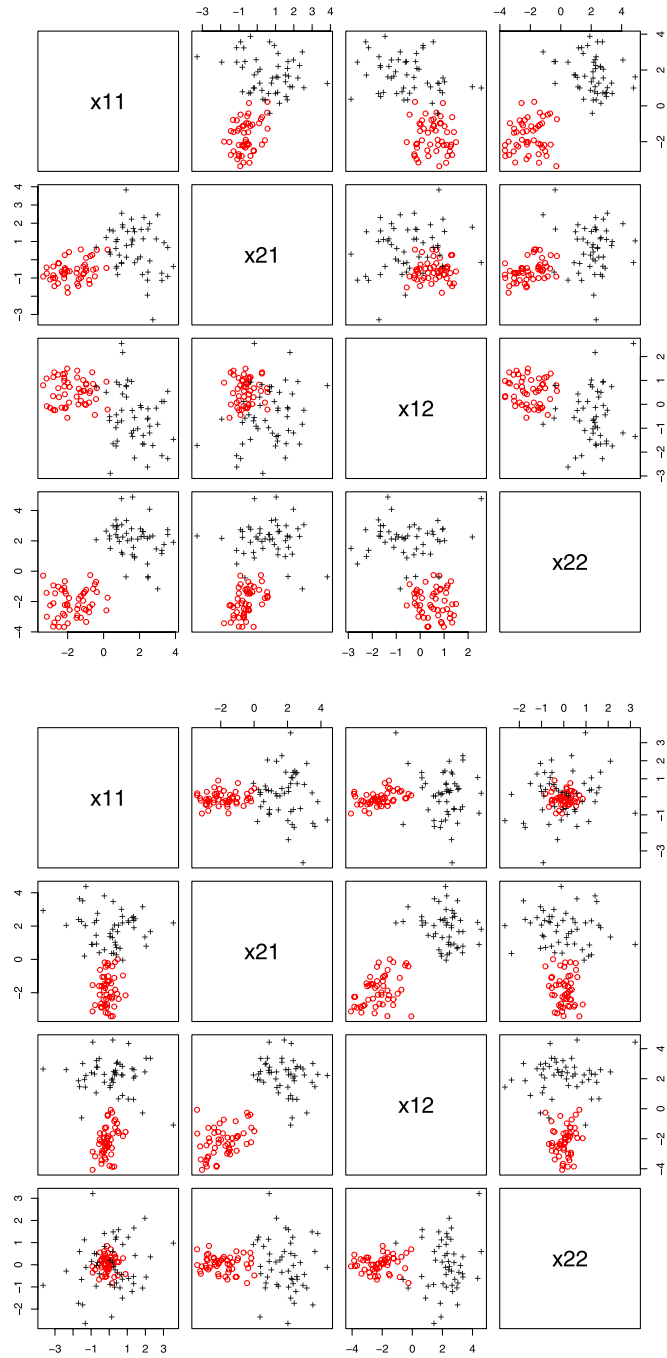

FIG. 2. *Scatterplot matrices for the reduced predictors estimated by folded-SIR (upper panel) and folded-DR (lower panel) for Example 1.*

Table 2
*Comparison between dimension-folding and conventional dimension reduction*

| Method | $n = 100$ | $n = 200$ | $n = 300$ | $n = 500$ | $n = 800$ |
|---|---|---|---|---|---|
| | $p_L = p_R = 5$ (benchmark distance = 2.586) | | | | |
| Folded-SIR | 1.057 | 0.716 | 0.580 | 0.432 | 0.343 |
| SIR | 1.865 | 1.806 | 1.758 | 1.759 | 1.753 |
| Folded-SAVE | 0.523 | 0.287 | 0.221 | 0.161 | 0.123 |
| SAVE | 1.615 | 1.294 | 1.089 | 0.757 | 0.579 |
| Folded-DR | 0.497 | 0.278 | 0.215 | 0.157 | 0.119 |
| DR | 1.596 | 1.289 | 1.075 | 0.747 | 0.574 |
| | $p_L = p_R = 10$ (benchmark distance = 2.586) | | | | |
| Folded-SIR | 1.924 | 1.246 | 0.984 | 0.750 | 0.577 |
| SIR | 2.626 | 2.142 | 2.051 | 1.963 | 1.921 |
| Folded-SAVE | 2.709 | 1.085 | 0.537 | 0.334 | 0.234 |
| SAVE | 2.753 | 2.677 | 1.956 | 1.605 | 1.406 |
| Folded-DR | 2.271 | 0.850 | 0.505 | 0.321 | 0.226 |
| DR | 2.753 | 2.503 | 1.871 | 1.593 | 1.392 |

**8. Application.** We now apply the dimension-folding methods to analyze the EEG data mentioned in the Introduction. The study involved two groups of subjects: an alcoholic group of 77 subjects and a control group of 45 subjects. Each subject was exposed to either one stimulus or two stimuli. During an exposure, the voltage values were measured from 64 channels of electrodes and for 256 time points (at 256 Hz per second). The 64 electrodes are placed at different locations on the subject's scalp. The stimuli were pictures chosen from a picture set. When two pictures were shown, they were displayed in either a matched condition, where two pictures were identical, or a unmatched condition, where they were different. Each subject had 120 trials under these three conditions: single stimulus, two matched stimuli and two unmatched stimuli. The primary interest was to study the association between alcoholism and the pattern of voltage values over times and channels.

To keep matters simple, in this paper, we use only part of the data set: we include only the single stimulus condition and, for each subject, we take the average of all the trials under that condition. That is, the portion of the data we use consists of $(\mathbf{X}_1, Y_1), \ldots, (\mathbf{X}_{122}, Y_{122})$, where $\mathbf{X}_i$ is a $256 \times 64$ matrix with each entry representing the mean voltage value of subject $i$ at a combination of a time point and a channel, averaged over all trials under the single stimulus condition, and $Y_i$ is a binary random variable indicating whether the $i$th subject is alcoholic ($Y_i = 1$) or nonalcoholic ($Y_i = 0$).

To apply the dimension-folding methods, we need to perform the spectral decomposition on the $p_L p_R \times p_L p_R = 16384 \times 16384$-dimensional matrix $\hat{\boldsymbol{\Sigma}}$, which is quite large. So, prior to dimension folding, we have implemented a somewhat heuristic pre-screening phase. Let $\mathbf{v}_1, \ldots, \mathbf{v}_{s_L}$ be the first $s_L$ eigenvectors of the matrix $E_n(\mathbf{X} - \bar{\mathbf{X}})(\mathbf{X} - \bar{\mathbf{X}})^T$ and $\mathbf{w}_1, \ldots, \mathbf{w}_{s_R}$ be the first $s_R$ eigenvectors of the matrix $E_n(\mathbf{X} - \bar{\mathbf{X}})^T(\mathbf{X} - \bar{\mathbf{X}})$. Let $\mathbf{V} = (\mathbf{v}_1, \ldots, \mathbf{v}_{s_L})$ and $\mathbf{W} = (\mathbf{w}_1, \ldots, \mathbf{w}_{s_R})^T$. Let $\mathbf{X}_i^* = \mathbf{V}^T \mathbf{X}_i \mathbf{W}$.

Using two sets of dimensions, $(s_L, s_R, d_L, d_R) = (30, 20, 2, 2)$ and $(30, 20, 1, 2)$, we apply folded-SIR, folded-SAVE and folded-DR to the pre-screened data set $(\mathbf{X}_1^*, Y_1), \ldots, (\mathbf{X}_n^*, Y_n)$. The results for $(d_L, d_R) = (2, 2)$ are presented in Figure 3, which contains two scatterplot matrices of the four predictors in the $2 \times 2$ matrices $\hat{\mathbf{a}}^T \mathbf{X}^* \hat{\mathbf{b}}$ obtained by folded-SAVE (left panel) and folded-DR (right panel). The four predictors are labeled as X11, X12, X21, X22 in the plots. A striking feature of these plots is that the EEG data for alcoholic cases (represented by red ○'s) show markedly less variation than those for nonalcoholic cases (represented by black +'s). This can be interpreted as indicating that the EEG patterns for the alcoholic subjects are more similar than those for the nonalcoholic cases. We also observe that folded-SAVE predictors show strong separation by variation, but no obvious separation by location, whereas folded-DR successfully separates the two clusters by both location and variation. The results for $(d_L, d_R) = (1, 2)$ are presented in Figure 4, which contains three scatterplots for the two predictors in the $1 \times 2$ matrix $\hat{\mathbf{a}}^T \mathbf{X}^* \hat{\mathbf{b}}$ obtained by folded-SIR (upper panel), folded-SAVE (lower-left panel) and folded-DR (lower-right panel). The two predictors are labeled as X11, X12 in the plots. From these plots, we observe the differences in performance of the three methods: folded-SIR works well in separating locations, folded-DR works well in separating variations, whereas folded-DR combines the advantages of both.

Of course, the ultimate purpose of dimension folding (or more generally, dimension reduction) is to assist regression or classification. Thus, the true test for the usefulness of our methods is whether they can help us to identify whether or not a person is alcoholic using his/her EEG data. For this reason, we have performed a classification analysis after dimension folding. For each $i = 1, \ldots, n$, we withhold the $i$th subject from the sample, treating it as the test set and the remaining 121 subjects as the training set. Based on each training set, we first carry out dimension folding (including pre-screening) and then apply quadratic discriminant analysis [Johnson and Wichern (2007), Chapter 11] to develop a classification rule. This classification rule is then used to classify the withheld subject. Using $(s_L, s_R, d_L, d_R) = (15, 15, 1, 2)$ and the folded-DR, we correctly predicted (as alcoholic or nonalcoholic) 97 out of the 122 cases; folded-SIR correctly classifies 94 out of 122 cases. We also compute the number of correct classifications using the conventional SIR, which gives 92 out 122 correct decisions. For the conventional SIR, we use $(s_L, s_R) = (9, 9)$ and $d = 1$. For all three methods, the ridge-regression-type inverse $(\hat{\boldsymbol{\Sigma}} + \varepsilon \mathbf{I}_{p_R p_L})^{-1}$ is used, with $\varepsilon = 0.5$. Noting that we have used only a portion of the data set, it is conceivable that an even stronger association could be established if the full data set were used.

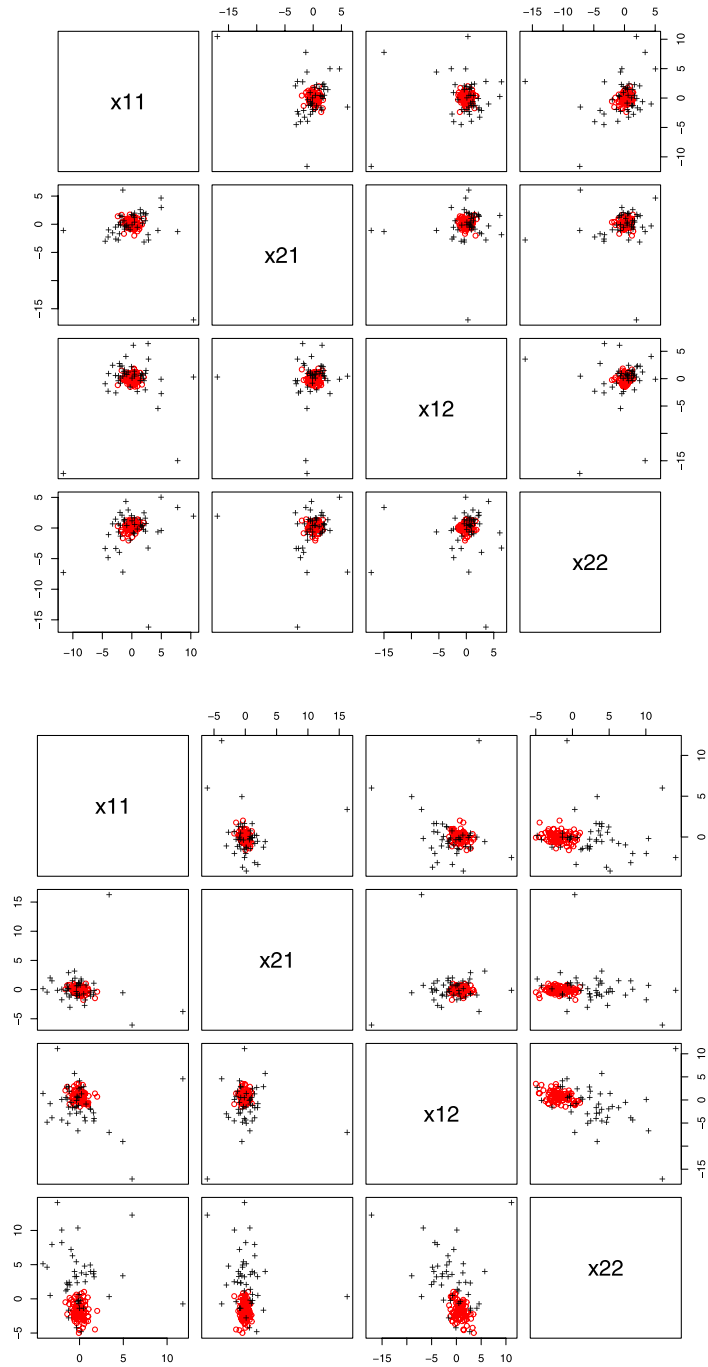

FIG. 3. *Scatterplot matrices for the four reduced predictors estimated by folded-SAVE (upper panel) and folded-DR (lower panel), for $(d_L, d_R) = (2, 2)$. Red $\circ$'s represent the alcoholic cases; black $+$'s represent the nonalcoholic cases.*

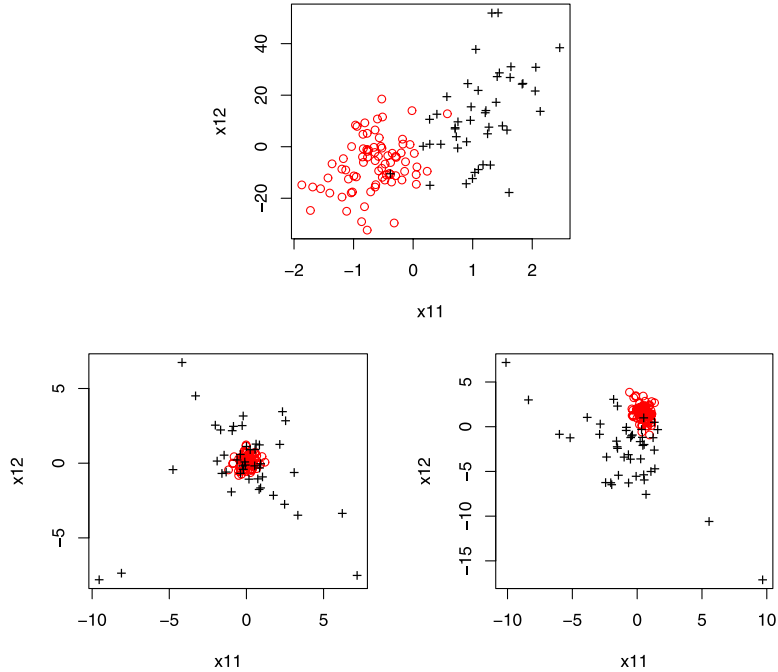

Fig. 4. *Scatterplots for the two reduced predictors estimated by folded-SIR (upper panel), folded-SAVE (lower-left panel) and folded-DR (lower-right panel), for $(d_L, d_R) = (1, 2)$. Red $\circ$'s represent the alcoholic cases; black $+$'s represent the nonalcoholic cases.*

The matrix structure preserved by dimension folding is helpful to gain further insights into the relationship between EEG patterns and alcoholism. In particular, the right dimension-folding subspace contains the weights of the channels that are important in predicting alcoholism, whereas the left dimension-folding subspace contains the principal patterns of how voltage varies over time in the important channels. These could provide important information for understanding how each part of the brain and the way it responds to stimuli are related to alcoholism.

**Acknowledgments.** We are thankful to three referees and an Associate Editor their exceptionally careful and insightful reviews that led to considerable improvement of this paper.

## REFERENCES

Billingsley, P. (1986). *Probability and Measure*, 2nd ed. Wiley, New York. MR0830424

BURA, E. and COOK, R. D. (2001). Estimating the structural dimension of regressions via parametric inverse regression. *J. R. Stat. Soc. Ser. B Stat. Methodol.* **63** 393–410. MR1841422

CHIAROMONTE, F. and MARTINELLI, J. (2002). Dimension reduction strategies for analyzing global gene expression data with a response. *Math. Biosci.* **176** 123–144. MR1869195

CHIAROMONTE, F. and COOK R. D. (2001). Sufficient dimension reduction and graphics in regression. *Ann. Inst. Statist. Math.* **54** 768–795. MR1954046

COOK, R. D. (1994). On the interpretation of regression plots. *J. Amer. Statist. Assoc.* **89** 177–189. MR1266295

COOK, R. D. (1996). Graphics for regressions with a binary response. *J. Amer. Statist. Assoc.* **91** 983–992. MR1424601

COOK, R. D. (1998). *Regression Graphics: Ideas for Studying Regressions through Graphics.* Wiley, New York. MR1645673

COOK, R. D. and LI, B. (2002). Dimension reduction for the conditional mean. *Ann. Statist.* **30** 455–474. MR1902895

COOK, R. D. and LI, B. (2004). Determining the dimension of iterative Hessian transformation. *Ann. Statist.* **32** 2501–2531. MR2153993

COOK, R. D., LI, B. and CHIAROMONTE, F. (2007). Dimension reduction without matrix inversion. *Biometrika* **94** 596–584. MR2410009

COOK, R. D., LI, B. and CHIAROMONTE, F. (2009). Envelope models for parsimonious and efficient multivariate linear regression (with discussion). *Statist. Sinica.* To appear.

COOK, R. D. and NI, L. (2005). Sufficient dimension reduction via inverse regression: A minimum discrepancy approach. *J. Amer. Statist. Assoc.* **100** 410–428. MR2160547

COOK, R. D. and WEISBERG, S. (1991). Discussion of "Sliced inverse regression for dimension reduction." *J. Amer. Statist. Assoc.* **86** 316–342. MR1137117

DONG, Y. and LI, B. (2009). Dimension reduction for nonelliptically distributed predictors: Second-order methods. *Biometrika.* Submitted. MR2509074

DUAN, N. and LI, K.-C. (1991). Slicing regression: A link-free regression method. *Ann. Statist.* **19** 505–530. MR1105834

JOHNSON, R. A. and WICHERN, D. W. (2007). *Applied Multivariate Statistical Analysis.* Pearson Prentice Hall, Upper Saddle River, NJ. MR2372475

FERRE, L. and YAO, A. F. (2005). Smooth function inverse regression. *Statist. Sinica* **15** 665–683. MR2233905

FUNG, K. F., HE, X., LIU, L. and SHI, P. (2002). Dimension reduction based on canonical correlation. *Statist. Sinica* **12** 1093–1113. MR1947065

HOERL, A. E. (1962). Application of ridge analysis to regression problems. *Chemical Engineering Progress* **58** 54–59.

LI, B. (2008). Comments on: Augmenting the bootstrap to analyze high dimensional genomic data. *Test* **17** 19–21. MR2393342

LI, B. and DONG, Y. (2009). Dimension reduction for nonelliptically distributed predictors. *Ann. Statist.* **37** 1272–1298. MR2509074

LI, B. and WANG, S. (2007). On directional regression for dimension reduction. *J. Amer. Statist. Assoc.* **102** 2143–2172. MR2354409

LI, B., WEN, S. and ZHU, L.-X. (2008). On a Projective Resampling method for dimension reduction with multivariate responses. *J. Amer. Statist. Assoc.* **103** 1177–1186. MR2462891

LI, B., ZHA, H. and CHIAROMONTE, C. (2005). Contour regression: A general approach to dimension reduction. *Ann. Statist.* **33** 1580–1616. MR2166556

LI, L. and LI, H. (2004). Dimension reduction methods for microarrays with application to censored survival data. *Bioinformatics* **20** 3406–3412.

LI, K.-C. (1991). Sliced inverse regression for dimension reduction (with discussion). *J. Amer. Statist. Assoc.* **86** 316–342. MR1137117

LI, K.-C. (1992). On principal Hessian directions for data visualization and dimension reduction: Another application of Stein's lemma. *J. Amer. Statist. Assoc.* **87** 1025–1039. MR1209564

MAGNUS, J. R. and NEUDECKER, H. (1979). The commutation matrix: Some properties and applications. *Ann. Statist.* **2** 381–394. MR0520247

MARQUARDT, D. W. (1970). Generalized inverses, ridge regression, biased linear estimation, and nonlinear estimation. *Technometrics* **12** 591–612.

TYEKUCHEVA, F. and CHIAROMONTE, F. (2008). Augmenting the bootstrap to analyze high-dimensional genomic data. *Test* 1–18. MR2393341

YE, Z. and WEISS, R. E. (2003). Using the bootstrap to select one of a new class of dimension reduction methods. *J. Amer. Statist. Assoc.* **98** 968–979. MR2041485

YIN, X. and COOK, R. D. (2002). Dimension reduction for the conditional $k$th moment in regression. *J. R. Stat. Soc. Ser. B Stat. Methodol.* **64** 159–175. MR1904698

YIN, X., LI, B. and COOK, R. D. (2008). Successive direction extraction for estimating the central subspace in a multiple-index regression. *J. Multivariate Anal.* **99** 1733–1757. MR2444817

ZHOU, J. (2009). Robust dimension reduction based on canonical correlation. *J. Multivariate Anal.* **100** 195–209. MR2460487

ZHU, L.-X. and FANG, K.-T. (1996). Asymptotics for kernel estimate of sliced inverse regression. *Ann. Statist.* **24** 1053–1068. MR1401836

DEPARTMENT OF STATISTICS
PENNSYLVANIA STATE UNIVERSITY
326 THOMAS BUILDING
UNIVERSITY PARK, PENNSYLVANIA 16802
USA
E-MAIL: bing@stat.psu.edu
         mzk136@psu.edu
         naomi@stat.psu.edu